\documentclass[10pt]{article}

\usepackage{graphicx}
\usepackage{amssymb}
\usepackage{epstopdf}
\usepackage{tikz}
\usepackage{lipsum}
\usepackage{amsfonts}
\usepackage{pgfplots}
\pgfplotsset{compat=1.18}

\definecolor{handwritten}{RGB}{140, 20, 80}

\usepackage{amsmath}
\usepackage[colorlinks=true,
            linkcolor=blue,
        citecolor=blue,
            urlcolor=blue]{hyperref}
\usepackage[nocompress]{cite}
\usepackage{amssymb}
\usepackage{comment}
\usepackage{color}
\usepackage{amsthm}

\def\phi{{\varphi}}

\DeclareSymbolFont{AMSb}{U}{msb}{m}{n}

\usepackage{indentfirst}
\usepackage{tikz}
\usepackage{pgfplots}
\pgfplotsset{compat=1.18}
\usepgfplotslibrary{groupplots}

\newtheorem{theorem}{Theorem}[section]
\newtheorem{proposition}[theorem]{Proposition}
\newtheorem{lemma}[theorem]{Lemma}
\newtheorem{remark}[theorem]{Remark}

\newtheorem{definition}[theorem]{Definition}

\newcommand\R{{\mathbb R}}

\title{Variational Thresholds and Multiplicity for a
Nonlocal $r$-Laplacian Logistic-Type Problem \footnote{This research was partially supported by CNPq grants from the Brazilian Federal Government and FAPEAM grants from the government of the State of Amazonas, Brazil.}}
\author{Thiago Cacau \thanks{Universidade Federal do Amazonas, Manaus-AM, Brazil, thiago.cacau@ufam.edu.br
 } \quad Alireza Khatib \thanks{Universidade Federal do Amazonas, Manaus-AM, Brazil, alireza@ufam.edu.br
 }  \quad
}
\date{}
\usetikzlibrary{arrows.meta,patterns}
\usepackage{float}
\begin{document}

 \maketitle
\begin{abstract}
This paper investigates the existence, multiplicity and asymptotic behavior of solutions to a class of nonlocal logistic-type problems with the \(r\)-Laplacian operator in the fast diffusion regime: %, with a real parameter $\lambda$:
\[
\begin{cases}
-\Delta_r u - \displaystyle\int_{\Omega} K(x,y)u(y)\,dy
= a(x)u +\lambda |u|^{q-2}u
\displaystyle\int_\Omega Q(x,y)|u(y)|^q\,dy,
& x\in \Omega, \\[0.4em]
u = 0, & x\in \partial\Omega.
\end{cases}
\]
The model has local diffusion, nonlocal dispersal, and nonlinear global interactions, and depends on a real parameter $\lambda$. For $\lambda \ge 0$, we establish the existence of nonnegative weak solutions via a variational principle. In the case \(\lambda<0\), we characterize two variational thresholds, \(-\infty<\widehat{\lambda}^*<\lambda^*<0\), which separate different solution intervals. For \(\lambda^*<\lambda<0\), the energy functional has a global minimizer with negative energy and a mountain-pass critical point with positive energy, yielding two distinct nontrivial weak solutions. The same holds at \(\lambda=\lambda^*\), where the global minimizer has zero energy. In a left neighborhood of \(\lambda^*\), two distinct nontrivial weak solutions persist, a local minimizer and a mountain-pass critical point, both with positive energy and the latter at a strictly higher level. For \(\lambda<\widehat{\lambda}^*\), only the trivial solution exists.
Finally, as $\lambda \to 0^-$, the norm of the global minimizers blows up. Along subsequences, selected nonnegative mountain-pass solutions converge strongly in $W_0^{1,r}(\Omega)$ as $\lambda \to 0$ to a nontrivial nonnegative weak solution of the problem at $\lambda = 0$. The associated mountain-pass levels converge to the same
limit from both sides of zero, and the corresponding
solution norms also converge to a common value.

\end{abstract}


%\tableofcontents

\section{Introduction}
Classical mathematical modeling in population dynamics developed from the logistic model proposed by Verhulst (1838), in which a term for the limited resources of the habitat regulates the growth rate of a population. This approach arose in contrast to the Malthusian model (1798), which assumed unlimited exponential growth \cite{verhulst1838}. As mathematical ecology advanced, models began to account for the spatial distribution and mobility of individuals within their habitat $\Omega$. Following Skellam \cite{skellam1951} and, later, Cantrell and Cosner \cite{cantrell2003}, short-distance local random movement came to be modeled by the classical Laplacian operator $-\Delta u$, which leads to local reaction-diffusion equations. Logistic reaction-diffusion problems involving nonlinear operators ($r \neq 2$), such as the $r$-Laplacian $\text{div}(\vert{}\nabla u\vert{}^{r-2}\nabla u)$, have been widely studied \cite{sabinalis2020, sabinadelis2022, oruganti2004}. In these population dynamics models, the exponent $r$ of the $r$-Laplacian operator sets the spatial dispersal regime of the species: slow diffusion for $r>2$ and fast diffusion for $1<r<2$ \cite{yang2022, wang2023, cong2016, rani2024, li2020}. In recent decades, nonlocal operators have been introduced where local diffusion may fall short: long-distance dispersal, jump-type movements (birth-jump processes), and competition for resources, which acts globally across the habitat \cite{alves2020, aydogmus2015, ma2018, ni2018, hillen2014, delgado2017}. 

We consider a nonlocal model combining local diffusion, nonlocal dispersal, and both cooperative and competitive nonlinear interactions, as studied in \cite{khatib2027}. The problem reads:
\begin{equation}
\begin{cases}\label{eq:problem}\tag{$P_\lambda$}
-\Delta_r u - \displaystyle\int_{\Omega}K(x,y)u(y)\,dy = a(x)u + \lambda |u|^{q-2}u \displaystyle \int_{\Omega}Q(x,y)|u(y)|^{q}\,dy & x \in \Omega \\ 
u=0, 
& x \in \partial \Omega,
\end{cases}
\end{equation}
where $\Omega \subset \R^N$ ($N \ge 1$) is a bounded domain with smooth boundary, $\lambda$ is a real parameter, and the $r$-Laplacian operator is defined by $\Delta_{r}u = \text{div}(\vert{}\nabla u\vert{}^{r-2}\nabla u)$. The exponents $r$ and $q$ present in problem \eqref{eq:problem} satisfy the following Sobolev restrictions 
\[ 
\max\left\{\frac{2N}{N+2},1\right\}<r<2 \le q < r^*, 
\]
where $r^*:= \frac{Nr}{N-r}$ denotes the critical Sobolev exponent when $r<N$ (and $r^*=+\infty$ if $N=1$). We assume the following hypotheses on the coefficients and interaction kernels:
\begin{itemize}
\item[($H_1$)] $a \in L^\infty(\Omega)$ satisfies $a(x) \ge a_0 > 0$ a.e. in $\Omega$;\item[($H_2$)] $K \in L^\infty(\Omega \times \Omega)$ is a symmetric, nonnegative kernel, that is, $K(y,x) = K(x,y) \ge 0$ a.e. in $\Omega \times \Omega$;
\item[($H_3$)] $Q \in L^\infty(\Omega \times \Omega)$ is a symmetric, strictly positive kernel, satisfying $Q(y,x)=Q(x,y) \ge q_0 > 0$ a.e. in $\Omega \times \Omega$.
\end{itemize}

The integral term 
\[
\int_\Omega K(x,y)u(y)\,dy
\]
describes a nonlocal dispersal process, where the kernel $K(x,y)$ represents the probability or migration rate of individuals from position $y$ to position $x$ \cite{cintra2025, sun2015, sun2021}. 

Depending on the sign of $\lambda \in \R$, the right-hand side of problem $\eqref{eq:problem}$ is cooperative for $\lambda>0$ and competitive, with logistic structure, for $\lambda<0$. The linear term $a(x)u$ models the intrinsic, spatially dependent reproduction rate of the population. The nonlocal nonlinear term 
\[
\lambda|u|^{q-2}u \int_{\Omega}Q(x,y)|u(y)|^q\,dy
\]
describes the global interactions among individuals throughout the habitat $\Omega$. When $\lambda<0$, the equation has a logistic character, and this nonlocal term represents global resource competition, limiting population growth based on the density everywhere in the domain. Conversely, when $\lambda > 0$, the model exhibits a cooperative regime, in which the nonlocal interaction term favors an increase in population density  \cite{alves2015, zhou2018, ma2020, garcia2009, shen2026}.

%Although the model studied in \cite{khatib2027} serves as motivation for the present work, our contribution lies in providing a variational analysis of problem \eqref{eq:problem} from a different perspective. This analysis highlights the structural role played by the parameter \(\lambda\) in the geometry of the associated energy functional, going beyond its interpretation as a mere constant. For $\lambda \ge 0$, the problem exhibits a Mountain Pass-type geometry. We employ the variational principle developed in \cite{moameni2020}, which enables us to obtain a nonnegative solution in a closed convex subset and subsequently establish that it is a weak solution in the whole space. The case $\lambda=0$ is treated as a particular instance of the cooperative regime, and the preservation of this variational geometry for $\lambda>0$ guarantees the existence of a corresponding nonnegative solution. In the competitive regime ($\lambda<0$), the use of the $r$-Laplacian operator enables a detailed study of the fibering map geometry. The presence of the $r$-Laplacian generates an algebraic interaction with the exponent $q$ of the global term that is structurally decisive. Such an interaction breaks the strict monotonicity of the energy functional, endowing the problem with the necessary geometry to harbor multiple solutions.

The model in \cite{khatib2027} motivates the present work, but we approach problem \eqref{eq:problem} differently, through the role the parameter $\lambda$ plays in the geometry of the energy functional. For this we use the variational structure of Silva \cite{Silva2020}, who established an abstract bifurcation result for equations involving homogeneous potential operators in reflexive Banach spaces. 

A direct application of Silva's result is not available,
since the quadratic form
\begin{equation}\label{quadratic-form}
\int_\Omega a(x)u^2\,dx
+\iint_{\Omega\times\Omega}
K(x,y)u(x)u(y)\,dx\,dy
\end{equation}
is positive on $E\setminus\{0\}$, but need not be positive
on sign-changing functions. Here
\[
X=W_0^{1,r}(\Omega),
\qquad
E:=\{u\in X:u\geq0\ \text{a.e. in }\Omega\}.
\]
Moreover, $E$ is a closed convex cone rather than a Banach subspace. We therefore adapt the variational arguments to this setting and verify the required properties directly for the nonlocal functional.

The results on existence, multiplicity, and nonexistence are stated below. 

\begin{theorem}\label{teo:1.1}
Assume $(H_1)-(H_3)$. Then, for every $\lambda \ge 0$, problem \eqref{eq:problem} admits a nontrivial nonnegative weak solution.
\end{theorem}
Let $I_\lambda$ denote the energy functional associated
with \eqref{eq:problem}, introduced in
Section~\ref{se.2}. The competitive regime is governed
by two variational thresholds,
\[
-\infty<\widehat{\lambda}^*<\lambda^*<0.
\]
The energy threshold $\lambda^*$ marks the onset of negative energy: for $\lambda<0$, the functional $I_\lambda$ takes negative values if and only if $\lambda>\lambda^*$. The Nehari threshold $\widehat{\lambda}^*$ is the lowest parameter for which the identity
\[
\langle I_\lambda'(u),u\rangle=0
\]
can hold for a nonzero function $u\in X$. Their precise variational definitions and attainment are established in Section \ref{se.4}; see Lemmas~\ref{lem:lambda*} and~\ref{lem:hat-lambda*}.
\begin{theorem}\label{teo:1.2}
Assume $(H_1)-(H_3)$. Then, for $\lambda \in (-\infty, 0)$, the following statements hold:
\begin{itemize}
    \item[\rm (i)] for each $\lambda \in (\lambda^*, 0)$, problem \eqref{eq:problem} admits two distinct weak solutions $v_\lambda, u_\lambda \in X$ satisfying $I_\lambda(u_\lambda) > 0$ and $\inf_{v\in X} I_\lambda(v) = I_\lambda(v_\lambda) < 0$;
    
    \item[\rm (ii)] for $\lambda = \lambda^*$, problem \eqref{eq:problem} admits two distinct nontrivial weak solutions $v_{\lambda^*}, u_{\lambda^*} \in X$ satisfying $I_{\lambda^*}(u_{\lambda^*}) > 0$ and $\inf_{v \in X} I_{\lambda^*}(v) = I_{\lambda^*}(v_{\lambda^*}) = 0$;

    \item[\rm (iii)] there exists $\varepsilon > 0$ such that, for each $\lambda \in (\lambda^* - \varepsilon, \lambda^*)$, problem \eqref{eq:problem} admits two distinct weak solutions $v_\lambda, u_\lambda \in X$ satisfying $I_\lambda(u_\lambda) >I_\lambda(v_\lambda) > 0$;

    %\item[\rm (iv)] if $\lambda < \lambda^*$, then $u \equiv 0$ is the unique
     % global minimizer of the energy functional $I_\lambda$ in $X$;
    
    \item[\rm (iv)] for $\lambda < \widehat{\lambda}^*$, problem \eqref{eq:problem} admits only the trivial solution.
\end{itemize}
\end{theorem}

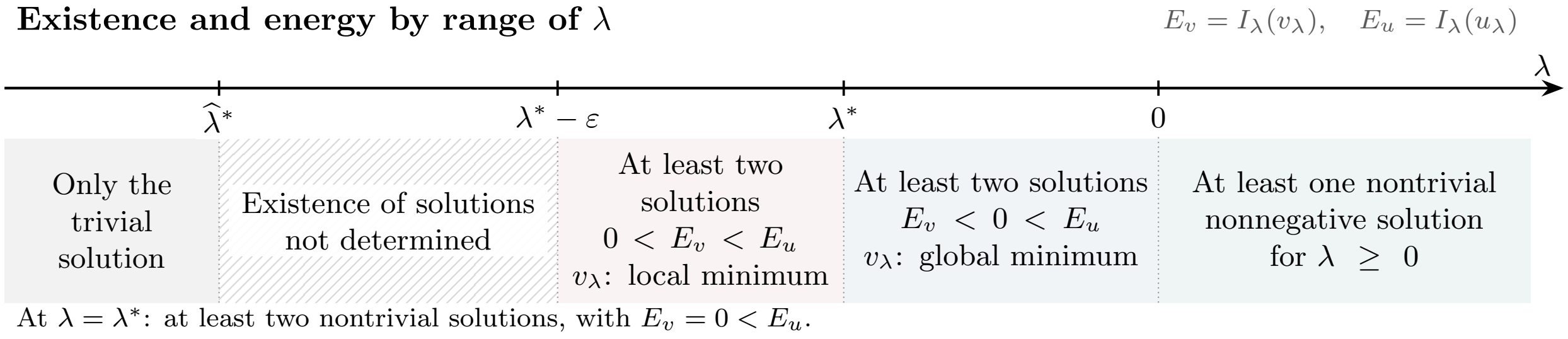
\begin{figure}[H]
\centering

\definecolor{mpred}{RGB}{164,48,63}
\definecolor{minblue}{RGB}{35,92,145}
\definecolor{coopteal}{RGB}{30,114,105}

\resizebox{\textwidth}{!}{%
\begin{tikzpicture}[
  x=1cm,
  y=1cm,
  font=\small,
  >=Stealth,
  regime/.style={
    align=center,
    inner sep=2pt,
    font=\small
  },
  ticklabel/.style={
    font=\small,
    fill=white,
    inner sep=2pt
  }
]
\node[anchor=west,font=\bfseries] at (0,10.35)
  {Existence and energy by range of $\lambda$};
\node[anchor=east,font=\footnotesize,text=black!65] at (16,10.35)
  {$E_v=I_\lambda(v_\lambda),\quad E_u=I_\lambda(u_\lambda)$};

\fill[black!5] (0,7.4) rectangle (2.25,9.12);
\fill[pattern=north east lines,pattern color=black!15]
  (2.25,7.4) rectangle (5.8,9.12);
\fill[mpred!6] (5.8,7.4) rectangle (8.8,9.12);
\fill[minblue!7] (8.8,7.4) rectangle (12.1,9.12);
\fill[coopteal!7] (12.1,7.4) rectangle (16,9.12);

\node[regime,text width=2.0cm] at (1.125,8.26)
  {Only the\\trivial solution};
\node[regime,text width=3.2cm,fill=white,fill opacity=.93,
      text opacity=1] at (4.025,8.26)
  {Existence of solutions\\not determined};
\node[regime,text width=2.8cm] at (7.3,8.26)
  {At least two solutions\\$0<E_v<E_u$\\$v_\lambda$: local minimum};
\node[regime,text width=3.1cm] at (10.45,8.26)
  {At least two solutions\\$E_v<0<E_u$\\$v_\lambda$: global minimum };
\node[regime,text width=3.6cm] at (14.05,8.26)
  {At least one nontrivial\\nonnegative solution\\for $\lambda\geq0$};

\draw[->,line width=.8pt] (0,9.65) -- (16.35,9.65)
  node[above left] {$\lambda$};
\foreach \xx in {2.25,5.8,8.8,12.1}{
  \draw[line width=.7pt] (\xx,9.57) -- (\xx,9.73);
  \draw[densely dotted,black!40] (\xx,7.4) -- (\xx,9.2);
}
\node[ticklabel] at (2.25,9.34) {$\widehat\lambda^*$};
\node[ticklabel] at (5.8,9.34) {$\lambda^*-\varepsilon$};
\node[ticklabel] at (8.8,9.34) {$\lambda^*$};
\node[ticklabel] at (12.1,9.34) {$0$};

\node[anchor=west,font=\footnotesize] at (0,7.22)
  {At $\lambda=\lambda^*$: at least two nontrivial solutions,
   with $E_v=0<E_u$.};
\end{tikzpicture}
}
\caption{Existence and multiplicity as $\lambda$ varies.}
\label{fig:parameter-regimes}
\end{figure}

Figure \ref{fig:parameter-regimes} shows the parameter ranges and energy properties from the main results.

We also study the asymptotic behavior of solutions as $\lambda\to0^{\pm}$, covering both the cooperative regime and the competitive regime. 

%\begin{theorem}\label{thm:asymptotic-behavior-zero}
%Assume $(H_1)-(H_3)$. Then the following assertions hold:
%\begin{enumerate}
%\item[(i)] For each $\lambda\in(\lambda^*,0)$, let
%$v_\lambda\in X$ be a global minimizer of $I_\lambda$
%such that $I_\lambda(v_\lambda)<0$. Then
%\[
%\|v_\lambda\|_X\to+\infty
%\qquad\text{as }\lambda\to0^-.
%\]
%Moreover, for $\lambda<0$ sufficiently close to zero, let
%$u_\lambda$ be the nonnegative mountain-pass solution
%chosen so that
%\[
%u_\lambda\in\mathcal M_\lambda\cap\mathcal N_\lambda^-,
%\]
%as provided by Lemma
%\ref{lem:bounded-mountain-pass-selection}.
%Then, for every sequence $\lambda_n\to0^-$,
%there exist a subsequence, still denoted by $(\lambda_n)$,
%and a function $u_0^-\in E\setminus\{0\}$ such that
%\[
%u_{\lambda_n}\to u_0^-
%\qquad\text{strongly in }X.
%\]
%Furthermore,
%\[
%I_0'(u_0^-)=0
%\qquad\text{and}\qquad
%I_0(u_0^-)=c_0.
%\]

%\item[(ii)] For each $\lambda>0$, let
%$u_\lambda\in E\setminus\{0\}$ be a nonnegative weak
%solution obtained through the mountain-pass construction
%of Theorem \ref{teo:1.1}, satisfying
%$I_\lambda(u_\lambda)=d_\lambda$. Then
%\[
%d_\lambda\to c_0
%\qquad\text{as }\lambda\to0^+.
%\]
%Moreover, for every sequence
%$(\lambda_n)\subset(0,\infty)$ with
%$\lambda_n\to0^+$, there exist a subsequence,
%still denoted by $(\lambda_n)$, and a function
%$u_0^+\in E\setminus\{0\}$ such that
%\[
%u_{\lambda_n}\to u_0^+
%\qquad\text{strongly in }X.
%\]
%Furthermore,
%\[
%I_0'(u_0^+)=0
%\qquad\text{and}\qquad
%I_0(u_0^+)=c_0.
%\]
%\end{enumerate}
%\end{theorem}

\begin{theorem}
\label{thm:asymptotic-behavior-zero}
Assume $(H_1)-(H_3)$. Then the following assertions hold.

\begin{enumerate}
\item[\rm (i)]
For each $\lambda\in(\lambda^*,0)$, let $v_\lambda$ be a global
minimizer of $I_\lambda$. Then
\[
\|v_\lambda\|_X\to+\infty
\qquad\text{as }\lambda\to0^-.
\]

\item[\rm (ii)]
For $\lambda<0$ sufficiently close to zero, the nonnegative mountain-pass solutions $u_\lambda$ can be chosen at level $c_\lambda$ so that the conclusions below hold. For $\lambda>0$, let $u_\lambda$ be a nonnegative mountain-pass solution at level $d_\lambda$ obtained in Theorem \ref{teo:1.1}. The corresponding levels satisfy
\[
\lim_{\lambda\to0^-}c_\lambda
=
\lim_{\lambda\to0^+}d_\lambda
=
c_0.
\]
For either one-sided limit, every sequence $\lambda_n\to0^\pm$ has a subsequence along which $u_{\lambda_n}$ converges strongly in $X$ to a nontrivial nonnegative weak solution $u_0$ of \eqref{eq:problem} at $\lambda=0$, with $I_0(u_0)=c_0$. In particular,
\[
\lim_{\lambda\to0^-}\|u_\lambda\|_{X}
=
\lim_{\lambda\to0^+}\|u_\lambda\|_{X}
=
\left(\frac{c_0}{1/r-1/2}\right)^{1/r}.
\]
\end{enumerate}
\end{theorem}

%The proof is in Section \ref{se.5}. Figure \ref{fig:asymptotics} show how the solution norms behave.
\begin{figure}[H]
\centering

\definecolor{mpred}{RGB}{151,49,64}
\definecolor{minblue}{RGB}{31,85,132}
\definecolor{axisgrey}{RGB}{77,82,91}

\resizebox{0.98\linewidth}{!}{%
\begin{tikzpicture}[
  x=1cm, y=1cm, font=\small, >=Stealth,
  mp/.style={circle, fill=mpred, inner sep=1.45pt},
  mini/.style={rectangle, draw=minblue, fill=white,
               inner sep=1.4pt, line width=.75pt},
  axis/.style={axisgrey, line width=.65pt},
  annotation/.style={font=\small, anchor=west, inner sep=0pt}
]

% Legend
\node[mp] at (1.05,5.16) {};
\node[annotation] at (1.25,5.16)
  {Mountain-pass solutions $u_\lambda$};
\node[mini] at (6.15,5.16) {};
\node[annotation] at (6.37,5.16)
  {Global minimizers $v_\lambda$};

% Axes and reference lines
\draw[axis,->] (0,0) -- (13.15,0)
  node[right,inner sep=2pt] {$\lambda$};
\draw[axis,->] (0,0) -- (0,4.73)
  node[above,inner sep=3pt] {$\|\cdot\|_X$};

\draw[axisgrey!50,densely dotted] (6.92,0) -- (6.92,4.67);
\draw[axisgrey!67,densely dashed,line width=.55pt]
  (0,1.58) -- (12.94,1.58);

\draw[axis] (-.08,1.58) -- (0,1.58);
\node[anchor=east,inner sep=3pt] at (-.08,1.58) {$s_0$};
\draw[axis] (6.92,-.065) -- (6.92,.065);
\node[below=3pt] at (6.92,0) {$0$};

% Schematic mountain-pass markers on both sides of zero
\foreach \px/\py in {
  1.10/2.40, 2.10/2.16, 3.10/2.22, 4.10/1.94,
  5.00/2.02, 5.78/1.79, 6.31/1.71, 6.60/1.63,
  7.24/1.53, 7.58/1.45, 8.16/1.36, 8.91/1.41,
  9.80/1.13, 10.80/1.20, 11.83/.96, 12.70/1.05
}{
  \node[mp] at (\px,\py) {};
}

% Schematic global-minimizer markers for negative parameters
\foreach \px/\py in {
  1.10/2.95, 2.10/3.12, 3.10/3.06, 4.10/3.46,
  5.00/3.39, 5.78/3.75, 6.31/4.09, 6.60/4.38
}{
  \node[mini] at (\px,\py) {};
}

\draw[minblue,line width=.75pt,->]
  (6.60,4.49) -- (6.60,4.77);
\fill (6.92,1.58) circle (1.7pt);

% Limit labels
\node[annotation,text=minblue] at (.70,4.16)
  {$\|v_\lambda\|_X\to+\infty$};
\node[annotation,text=mpred] at (8.28,2.41)
  {$\|u_\lambda\|_X\to s_0$ from both sides};
\node[annotation,text=axisgrey] at (8.28,3.60)
  {$\displaystyle s_0=
    \left(\frac{c_0}{1/r-1/2}\right)^{1/r}$};

\draw[axis,->] (1.25,-.57) -- (6.54,-.57);
\node[fill=white,inner sep=4pt] at (3.90,-.57)
  {$\lambda\to0^-$};
\draw[axis,<-] (7.30,-.57) -- (12.59,-.57);
\node[fill=white,inner sep=4pt] at (9.95,-.57)
  {$\lambda\to0^+$};

\node[text=axisgrey,font=\footnotesize] at (6.92,-1.14)
  {The minimax levels satisfy $c_\lambda\to c_0$ from the left
   and $d_\lambda\to c_0$ from the right.};

\end{tikzpicture}%
}

\caption{Schematic behavior of the solution norms as $\lambda\to0$. The norms of the global minimizers diverge as $\lambda\to0^-$, while the mountain-pass solutions have the same limiting norm from both sides.}
\label{fig:asymptotics}
\end{figure}
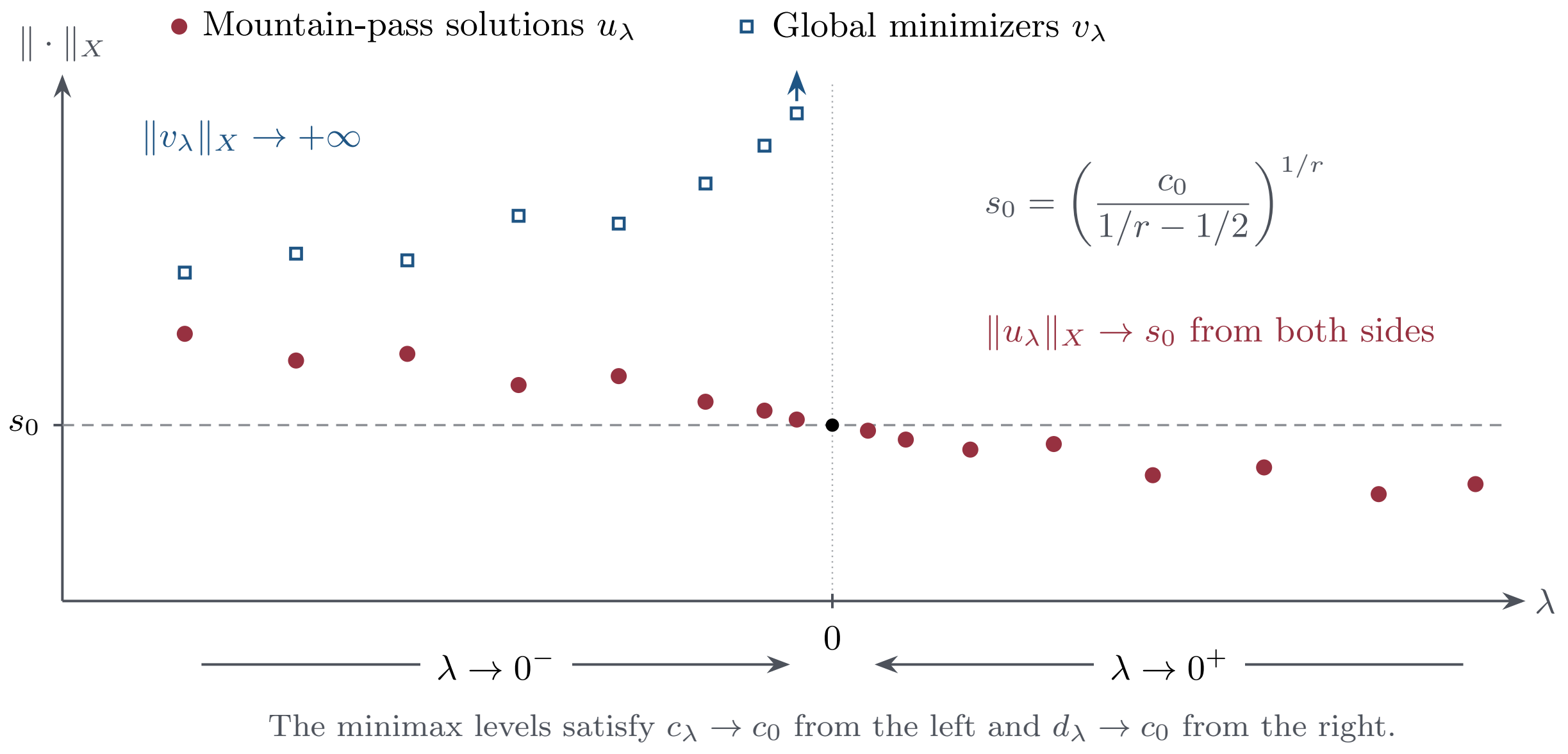

\medskip\noindent\textbf{Main contributions.} 
First, for $\lambda\geq0$, we construct a nontrivial
nonnegative weak solution through a convex restriction
argument and an auxiliary Dirichlet problem.
Second, for $\lambda<0$, we identify two variational
thresholds and establish multiplicity and nonexistence
results, including the persistence of two positive-energy
solutions immediately below the energy threshold.
Third, we describe the behavior as $\lambda\to0$:
global minimizers become unbounded as $\lambda\to0^-$,
whereas selected nonnegative mountain-pass solutions
converge along subsequences to solutions at $\lambda=0$.
Their energy levels and norms have common limits
from both sides.

%First, for $\lambda\ge0$, we use the variational principle proposed by \cite{moameni2020} to obtain a nonnegative weak solution in the closed convex subset $E$, and then show that it is a weak solution in the whole space $X$. Second, for $\lambda<0$, we characterize the two extremal parameters and obtain multiplicity and nonexistence results. Third, we study the behavior of the solutions as $\lambda\to0$. For selected nonnegative mountain-pass solutions, the minimax levels converge to the mountain-pass level $c_0$ of $I_0$, and along subsequences the solutions converge strongly in $X$ to a nontrivial nonnegative solution $u_0$ of the problem with $\lambda=0$ and $I_0(u_0)=c_0$, and the norms converge, from both sides and along the whole family, to $\big(c_0/\left(\tfrac1r-\tfrac12\right)\big)^{1/r}$.  In contrast, the $X$-norms of the global minimizers blow up as $\lambda \to 0^-$. We do not claim that the limit is unique, nor that the limits from the two sides of $\lambda=0$ coincide.

%{\color{red}
%The paper is organized as follows. In section \ref{se.2}, we fixed notations and results preliminary. In section \ref{se.3}, we used a variational principle to prove that existence nonnegative weak solution in $X$ to problem \ref{eq:problem} in case $\lambda \ge 0$ regime. In section \ref{se.4}, we treat the case $\lambda<0$ regime, where we identification two parameter extremal to prove the existence multiple of solution for different intervals of parameter. Finally, in \ref{se.5}, we show behavior asymptotic as $\lambda \to 0$.
%}

\section{Preliminaries}\label{se.2}
This section fixes the notation and collects preliminary results used throughout the paper.

%\begin{theorem}[Variational Principle]\label{principle}
%Let $\Psi : X \to \mathbb{R} \cup \{+\infty\}$ be a proper convex and lower semi-continuous function and, let $E$ be a closed and convex subset of $X$. Assume that $\Psi$ is G\^ateaux differentiable at each point of $E$ and $\Phi \in C^1(X; \mathbb{R})$. If the following two assertions hold
%\begin{enumerate}
%    \item[(i)] The functional $I_E : X \to \mathbb{R} \cup \{+\infty\}$ defined as $I_{E,\lambda}(u) = \Psi_E(u) - \Phi(u)$ 
%    has a critical point at $\bar u$;
%    \item[(ii)] There exists $\bar v \in E$ such that $\Psi'(\bar v) =\Phi'(\bar u)$.
%\end{enumerate}
%Then $\bar u$ is a solution of the problem 
%\[
%\Psi'(\bar u) = \Phi'(\bar u).
%\]
%\end{theorem}

Recall that \(X=W_0^{1,r}(\Omega)\). We endow it with the norm
\[
\|u\|_X:=\left(\int_\Omega|\nabla u|^r\,dx\right)^{1/r},
\]
which is equivalent to the usual norm of \(W_0^{1,r}(\Omega)\) by the Poincaré inequality. Assume that 
\[
\max\left\{\frac{2N}{N+2},1\right\} < r < 2 \le q.
\]
The energy functional $I_\lambda : X \to \mathbb{R}$ associated with problem \eqref{eq:problem} is given by
\[
I_{\lambda}(u)=\Psi(u) -\Phi_{\lambda}(u),
\]
where the functionals $\Psi, \Phi_\lambda: X \to \R$ are defined by 
\[
\Psi(u)= \frac{1}{r}\int_{\Omega} |\nabla u|^rdx, \quad \Phi_\lambda(u)=\frac{1}{2}B(u)+\frac{\lambda}{2q}C(u). 
\]
with 
\[
B(u)=\int_\Omega a(x)u^2\,dx +\iint_{\Omega\times\Omega}K(x,y)u(x)u(y)\,dx\,dy, \quad C(u):=\iint_{\Omega\times\Omega} Q(x,y)|u(x)|^q|u(y)|^q\,dx\,dy.
\]
Let $E$ be the closed convex cone in $X$. We introduce the functional $\Psi_{E} : X \to \mathbb{R} \cup \{+\infty\}$ given by
\[
\Psi_{E}(u) := 
\begin{cases} 
\Psi(u), & \text{if } u \in E, \\ 
+\infty, & \text{if } u \notin E. 
\end{cases}
\]
The functional $I_{E,\lambda} : X \to \mathbb{R} \cup \{+\infty\}$ is the energy functional restricted to $E$, defined by
\[
I_{E,\lambda}(u) := \Psi_{E}(u) - \Phi_\lambda(u).
\]

The next two lemmas are used repeatedly throughout the paper.

\begin{lemma}\label{lem:funcional-c1}
Assume $(H_1)-(H_3)$. Then $\Phi_\lambda : X \to \mathbb{R}$ is well defined and of class $C^1$, and its derivative is given by
\[
\langle \Phi_\lambda'(u), v \rangle = \int_\Omega a(x)uv \, dx + \iint_{\Omega \times \Omega} K(x,y) u(y) v(x) \, dx \, dy + \lambda \int_\Omega |u(x)|^{q-2} u(x) v(x) \left( \int_\Omega Q(x,y) |u(y)|^q \, dy \right) dx
\]
for all $u, v \in X$. Consequently, the energy functional $I_\lambda = \Psi - \Phi_\lambda$ is also well defined and belongs to $C^1(X, \mathbb{R})$.
\end{lemma}

\begin{proof} We first show that $\Phi_\lambda(u)$ is well defined. Since $a \in L^\infty(\Omega)$ and $K, Q \in L^\infty(\Omega \times \Omega)$, Hölder's inequality gives
\[
|\Phi_\lambda(u)| \le \frac{1}{2}\|a\|_\infty \|u\|_{2}^2 + \frac{1}{2}\|K\|_\infty \|u\|_{1}^2 + \frac{|\lambda|}{2q}\|Q\|_\infty \|u\|_{q}^{2q}.
\]
Since $r^* > q \ge 2$, the continuous Sobolev embeddings give $\|u\|_{s} \le C_s \|u\|_X$ for all $1 \le s < r^*$, so $|\Phi_\lambda(u)|  < \infty$. The expression for $\langle \Phi_\lambda'(u), v \rangle$ stated in the lemma follows from a standard computation of the G\^ateaux derivative, using the symmetry of the kernels $K$ and $Q$.

It remains to show that $\Phi_\lambda' : X \to X^*$ is continuous. Let $(u_n) \subset X$ be a sequence converging strongly to $u$ in $X$, and let $v \in X$ be arbitrary. By continuous Sobolev embeddings, $u_n \to u$ strongly in $L^s(\Omega)$ for all $1 \le s < r^*$. In particular, the convergence $u_n \to u$ in $L^q(\Omega)$ implies that $|u_n|^q \to |u|^q$ in $L^1(\Omega)$ and that $|u_n|^{q-2}u_n \to |u|^{q-2}u$ strongly in $L^{\frac{q}{q-1}}(\Omega)$. 

For the terms involving $a$ and $K$, Hölder's inequality and Fubini's theorem give
\begin{align*}
\left| \int_\Omega a(x) (u_n - u)v \, dx \right| &\le \|a\|_\infty \|u_n - u\|_2 \|v\|_2 \to 0, \\
\left| \iint_{\Omega\times\Omega} K(x,y) (u_n(y) - u(y)) v(x) \, dx \, dy \right| &\le \|K\|_\infty \|v\|_1 \|u_n - u\|_1 \to 0.
\end{align*}

For the nonlocal term involving $Q$, we define 
\[
w_n(x) := \int_{\Omega} Q(x,y)|u_n(y)|^q \, dy \quad \text{and} \quad w(x) := \int_{\Omega} Q(x,y)|u(y)|^q \, dy.
\]
Then
\[
\|w_n - w\|_\infty \le \|Q\|_\infty \big\| |u_n|^q - |u|^q \big\|_1 \to 0.
\]
Adding and subtracting $|u|^{q-2}u v w_n$, we get 
\begin{align*}
\int_\Omega \Big| |u_n|^{q-2}u_n w_n - |u|^{q-2}u w \Big| |v| \, dx \le \|w_n\|_\infty \big\| |u_n|^{q-2}u_n - |u|^{q-2}u \big\|_{\frac{q}{q-1}} \|v\|_q + \|w_n - w\|_\infty \|u\|_q^{q-1} \|v\|_q.
\end{align*}
Since $\|w_n\|_\infty$ is uniformly bounded and both differences on the right-hand side tend to zero, the integral tends to zero as $n \to \infty$. Hence $\langle \Phi_\lambda'(u_n), v \rangle \to \langle \Phi_\lambda'(u), v \rangle$.

The same estimates bound $|\langle \Phi_\lambda'(u_n) - \Phi_\lambda'(u), v \rangle|$ by a sum of terms that tend to zero, each multiplied by the norm of $v$ in $L^1(\Omega)$, $L^2(\Omega)$, and $L^q(\Omega)$. Since $\Vert{}v\Vert{}_{s} \le C_s \Vert{}v\Vert{}_X$ for all $1 \le s < r^*$, these estimates are uniform on bounded sets in $X$. Taking the supremum over $\Vert{}v\Vert{}_X \le 1$, we obtain
\[
\|\Phi_\lambda'(u_n) - \Phi_\lambda'(u)\|_{X^*} \to 0.
\]
Since the G\^ateaux derivative is continuous, it is also the Fréchet derivative, so $\Phi_\lambda \in C^1(X, \mathbb{R})$. Finally, since the $r$-Laplacian part $\Psi(u) = \frac{1}{r}\|u\|_X^r$ is a standard functional known to be of class $C^1(X, \mathbb{R})$, we conclude that $I_\lambda = \Psi - \Phi_\lambda \in C^1(X, \mathbb{R})$.

\end{proof}

%To simplify notation, we define
%\[
%\begin{aligned}
%B(u)
%&:=\int_\Omega a(x)u^2\,dx
%+\iint_{\Omega\times\Omega}K(x,y)u(x)u(y)\,dx\,dy,\\
%C(u)
%&:=\iint_{\Omega\times\Omega}
%Q(x,y)|u(x)|^q|u(y)|^q\,dx\,dy.
%\end{aligned}
%\]

\begin{lemma}\label{lem:convergences} 
Assume $(H_1)-(H_3)$. If $u_n \rightharpoonup u$ weakly in $X$, then
\begin{gather}
     B(u_n) \to B(u), \,\,C(u_n) \to C(u); \\ 
    \langle \Phi_\lambda'(u_n), u_n - u \rangle \to 0.
\end{gather}
\end{lemma}

\begin{proof}Since $u_n \rightharpoonup u$ weakly in $X$, the sequence $(u_n)$ is bounded. By compact Sobolev embeddings, $u_n \to u$ strongly in $L^s(\Omega)$ for all $1 \le s < r^*$. Consequently, we have $u_n^2 \to u^2$ and $|u_n|^q \to |u|^q$ strongly in $L^1(\Omega)$. By adding and subtracting appropriate terms and applying Hölder's inequality, we obtain
\begin{align*}
|B(u_n) - B(u)| &\le C_1 \Big( \|u_n^2 - u^2\|_1 + \big(\|u_n\|_1 + \|u\|_1\big)\|u_n - u\|_1 \Big), \\
|C(u_n) - C(u)| &\le C_2 \big( \|u_n\|_q^q + \|u\|_q^q \big) \big\| |u_n|^q - |u|^q \big\|_1,
\end{align*}
where $C_1 > 0$ depends only on $\|a\|_\infty$ and $\|K\|_\infty$, while $C_2 > 0$ depends only on $\|Q\|_\infty$. Since the sequence $(u_n)$ is bounded in both $L^1(\Omega)$ and $L^q(\Omega)$, and all difference norms on the right-hand side vanish as $n \to \infty$, we conclude that $B(u_n) \to B(u)$ and $C(u_n) \to C(u)$. 
Testing the derivative against $u_n-u$ and applying Hölder's inequality and Fubini's theorem to each term, we obtain 
\begin{align*}
|\langle \Phi_\lambda'(u_n), u_n - u \rangle| &\le \|a\|_\infty \|u_n\|_2 \|u_n - u\|_2 + \|K\|_\infty \|u_n\|_1 \|u_n - u\|_1 + |\lambda|\|Q\|_\infty \|u_n\|_q^{2q-1} \|u_n - u\|_q \\
&\le \overline{C} \Big(  \|u_n\|_2\|u_n - u\|_2 +  \|u_n\|_1\|u_n - u\|_1 +  \|u_n\|_q^{2q-1}\|u_n - u\|_q \Big).
\end{align*}
Since the sequence $(u_n)$ is uniformly bounded in $L^s(\Omega)$ and $u_n \to u$ strongly in these spaces, the right-hand side vanishes as $n \to \infty$. This completes the proof.

\end{proof}

\begin{remark}\label{rem:IE_not_C1}
The extended functional $I_{E,\lambda}$ is not of class $C^1(X, \mathbb{R})$, because of the indicator function of the closed convex set, which takes the value $+\infty$ outside this set. This nonsmoothness rules out classical critical point theory and is why we use a nonsmooth variational framework.
\end{remark}

%De agora em diante, ao longo de todo o restante deste artigo, assumiremos que os núcleos $K$ e $Q$ são simétricos q.t.p. em $\Omega \times \Omega$, isto é,
%\[
%K(x,y) = K(y,x) \quad \text{e} \quad Q(x,y) = Q(y,x) \qquad \text{q.t.p. em } \Omega \times \Omega.
%\]

\section{The Nonnegative Case: $\lambda \ge 0$}\label{se.3}
This section treats the cooperative regime, \(\lambda>0\), and the limiting case \(\lambda=0\). Instead of the technical arguments commonly used to establish the positivity of solutions, we apply the variational principle of \cite{moameni2020} directly. We first recall some preliminary concepts from \cite{szulkin1986}, which are used throughout the paper, and then state a version of the principle adapted to our problem.

\begin{definition}[Critical Point]\label{critical_point}
Let $X$ be a real Banach space, let $\Phi \in C^1(X, \mathbb{R})$, and let $\Psi : X \to (-\infty, +\infty]$ be proper \rm (i.e., $\operatorname{Dom}(\Psi) \neq \emptyset$), convex, and lower semicontinuous. Consider the functional
\begin{equation}\label{H} \tag{H}
I := \Psi - \Phi.
\end{equation}
A point $u \in X$ is said to be a \textbf{critical point} of $I$ if $u \in \operatorname{Dom}(\Psi)$ and
\begin{equation*}
\Psi(v) - \Psi(u) \ge \langle\Phi'(u), v-u \rangle \quad \forall v \in X.
\end{equation*}
\end{definition}

\begin{definition}[Palais-Smale Condition]\label{PS-condition}
Let $X$ a Banach space. We say that $I$ of the form \eqref{H} satisfies the \textbf{Palais-Smale compactness condition} \rm(PS) if every sequence $(u_n) \subset X$ such that
\begin{itemize}
    \item[(i)]  $I(u_n) \to c$,
    \item[(ii)]  $\Psi(v) - \Psi(u_n) \ge  \langle\Phi'(u_n), v-u_n \rangle  -\epsilon_n \|v - u_n\|_X, \quad \forall v \in X$,
\end{itemize}
where $\epsilon_n \to 0$, has a convergent subsequence.
\end{definition}

\begin{theorem}[Nonsmooth Mountain-Pass Theorem]\label{Mountain_Pass}
Let $X$ a Banach space, and let $I : X \to (-\infty, +\infty]$ be a functional of the form \rm (\ref{H}) satisfying the Palais-Smale condition. Assume that 
\begin{itemize}
    \item[(i)]  $I(0) = 0$;
    \item[(ii)]  there exists $e \in X$ such that $I(e) \le 0$;
    \item[(iii)]  there exists some $\rho$ such that $0 < \rho < \|e\|_X$ and for every $u \in X$ with $\|u\|_X = \rho$ one has $I(u) > 0$.
\end{itemize}
Then $I$ has a critical value $c$ given by
\[
c = \inf_{\gamma \in \Gamma} \sup_{t \in [0,1]} I(\gamma(t)),
\]
where $\Gamma = \left\{\gamma \in C([0, 1], X) : \gamma(0) = 0, \; \gamma(1) = e\right\}$.
\end{theorem}

\begin{theorem}\label{adapted_principle}
Let $E$ be a nonempty, closed and convex subset of $X$. Assume that
\begin{enumerate}
    \item[(i)] The functional $I_{E,\lambda} : X \to \R \cup \{+\infty\}$ has a critical point $\bar u \in X$ in the sense of Definition \ref{critical_point};
    \item[(ii)] There exists $\bar{v} \in E$, with $\bar v = 0$ on $\partial \Omega$, such that
    \begin{equation}\label{itemii-principle}
    -\Delta_r \bar{v} = a(x)\bar{u} + \int_{\Omega} K(x,y)\bar{u}(y)\,dy + \lambda |\bar{u}|^{q-2}\bar{u} \int_{\Omega} Q(x,y)|\bar{u}(y)|^q\,dy 
    \end{equation}
    in the weak sense.
\end{enumerate}
Then $\bar{u} \in E$ is a weak solution of problem \eqref{eq:problem}.
\end{theorem}

\begin{proof}
Let $g(\bar u)$ denote the right-hand side in \eqref{itemii-principle},
\[
g(\bar u)=a(x)\bar{u} + \int_{\Omega} K(x, y)\bar{u}(y)\,dy + \lambda |\bar{u}|^{q-2}\bar{u} \int_{\Omega} Q(x, y)|\bar{u}(y)|^q\,dy.
\]
By Lemma \ref{lem:funcional-c1}, $\Phi_\lambda \in C^1(X,\R)$ and
\[
\langle \Phi_\lambda'(\bar u), w \rangle = \int_\Omega g(\bar u)\, w \, dx \quad \text{for all } w \in X.
\]
By (i), $\bar u$ is a critical point of $I_{E,\lambda}$, so $\bar u\in E$ and, by Definition \ref{critical_point},
\begin{equation}\label{4.4-1}
\Psi(v) - \Psi(\bar u) \ge \langle\Phi_\lambda'(\bar u), v - \bar u \rangle \quad \forall v \in E.
\end{equation}
Taking $v = \bar{v}$ gives
\begin{equation}\label{4.1-2}
\frac{1}{r}\int_{\Omega} |\nabla \bar{v}|^r\,dx - \frac{1}{r}\int_{\Omega} |\nabla \bar{u}|^r\,dx \ge \int_{\Omega} g(\bar{u})(\bar{v} - \bar{u})\,dx.
\end{equation}
The weak formulation in (ii) holds for every test function in $X$, so we may choose $\eta = \bar{v} - \bar{u}$, which gives
\begin{equation}\label{4.1-3}
\int_{\Omega} |\nabla \bar{v}|^{r-2}\nabla \bar v \cdot \nabla (\bar{v} - \bar{u})\,dx = \int_{\Omega} g(\bar{u})(\bar{v} - \bar{u})\,dx.
\end{equation}
Substituting \eqref{4.1-3} into \eqref{4.1-2}, we get
\begin{equation}\label{4.1-4}
\frac{1}{r}\left(\int_\Omega |\nabla \bar v|^rdx-\int_\Omega |\nabla \bar u|^rdx\right)
\ge
\int_\Omega
|\nabla\bar v|^{r-2}\nabla\bar v \cdot
(\nabla\bar v-\nabla\bar u)\,dx.
\end{equation}
Since $\Psi$ is convex and differentiable, we also have
\begin{equation*}
\frac{1}{r}\left(\int_\Omega |\nabla \bar u|^rdx-\int_\Omega |\nabla \bar v|^rdx \right)
\ge
\int_\Omega
|\nabla\bar v|^{r-2}\nabla\bar v\cdot
(\nabla\bar u-\nabla\bar v)\,dx. 
\end{equation*}
Thus equality holds in \eqref{4.1-4}, and the equality case of the convexity inequality yields
\[
\nabla\bar u=\nabla\bar v
\quad\text{a.e. in }\Omega.
\]
Because $\bar u-\bar v\in W_0^{1,r}(\Omega)$, it follows that 
\[
\bar u=\bar v \quad\text{a.e. in }\Omega.
\] 
Therefore, $\bar u\in E$ is a weak solution of problem \eqref{eq:problem}.

\end{proof}

\begin{lemma}\label{lem:PS-IE}
Assume $(H_1)-(H_3)$. For every fixed $\lambda\geq 0$ and every $d\in\mathbb{R}$, the functional $I_{E,\lambda}$ satisfies the Palais-Smale condition at level $d$.
\end{lemma}

\begin{proof}
Let $(u_n) \subset E$ be a Palais-Smale sequence at level $d$, that is, 
\[
I_{E,\lambda}(u_n) \to d,
\]
\begin{equation}\label{4.6-1}
\Psi(v)-\Psi(u_n) \ge \langle\Phi'_\lambda(u_n), v - u_n \rangle - \varepsilon_n \|v - u_n\|_X, \quad \text{ for all } v \in E,
\end{equation}
where $\varepsilon_n \to 0^+$. For $n$ sufficiently large, we have
\begin{equation}\label{4.6-2}
I_{E,\lambda}(u_n) \le |d| + 1.
\end{equation}
Let $t > 0$. Since $u_n \in E$, the function $v_t = (1+t)u_n$ satisfies $v_t \ge 0$, so $v_t \in E$. Taking $v = v_t$ in \eqref{4.6-1} and dividing by $t > 0$, we obtain
\[
\frac{(1+t)^r - 1}{rt}\|u_n\|_X^r \ge \langle\Phi'_\lambda(u_n), u_n \rangle - \varepsilon_n \|u_n\|_X.
\]
Since $\big((1+t)^r - 1\big)/t \to r$ as $t \to 0^+$, this gives
\begin{equation}\label{4.6-3}
\|u_n\|_X^r \ge \langle\Phi'_\lambda(u_n), u_n \rangle - \varepsilon_n \|u_n\|_X.
\end{equation}
Note that
%\[
%\langle\Phi'_\lambda(u_n), u_n \rangle = \int_\Omega a(x)u_n^2\,dx + \iint_{\Omega\times\Omega} K(x,y)u_n(x)u_n(y)\,dx\,dy + \lambda \iint_{\Omega\times\Omega} Q(x,y)|u_n(x)|^q|u_n(y)|^q\,dx\,dy.
%\]
\[
\langle\Phi'_\lambda(u_n), u_n \rangle = B(u_n) + \lambda C(u_n).
\]
Multiplying inequality \eqref{4.6-3} by $-1/2$ and adding the results to \eqref{4.6-2}, we get
%\[
%\left(\frac{1}{r} - \frac{1}{2}\right) \|u_n\|_X^r + \lambda \left(\frac{1}{2} - %\frac{1}{2q}\right) \iint_{\Omega\times\Omega} Q(x,y)|u_n(x)|^q|u_n(y)|%^q\,dx\,dy \le d + 1 + \frac{1}{2}\varepsilon_n \|u_n\|_X.
%\]
\[
\left(\frac{1}{r} - \frac{1}{2}\right) \|u_n\|_X^r + \lambda \left(\frac{1}{2} - \frac{1}{2q}\right) C(u_n) \le |d| + 1 + \frac{1}{2}\varepsilon_n \|u_n\|_X.
\]
Since $\varepsilon_n\to0$, we may assume $\varepsilon_n\leq1$ for all large $n$. Using $\lambda\geq0$, $q\geq2$, and $C(u_n)\geq0$, we deduce
\[
\left(\frac{1}{r} - \frac{1}{2}\right) \|u_n\|_X^r \le |d| + 1 + \frac{1}{2}\|u_n\|_X.
\]
Since $1<r<2$, this inequality shows that $(u_n)$ is bounded in $X$.

Since $X$ is reflexive and $(u_n)$ is bounded in $X$, there exists $u \in X$ such that, up to a subsequence,
\[
u_n \rightharpoonup u \quad \text{weakly in } X.
\]
The set $E$ is closed and convex, hence weakly closed, so $u \in E$. Taking $v=u\in E$ in \eqref{4.6-1} gives
\[
\frac{1}{r}\|u_n\|_X^r \le \frac{1}{r}\|u\|_X^r + \langle\Phi'_\lambda(u_n), u_n - u \rangle + \varepsilon_n \|u_n - u\|_X.
\]
By Lemma \ref{lem:convergences}, $\langle \Phi_\lambda'(u_n),u_n-u\rangle\to0$. Since $(u_n)$ is bounded and $\varepsilon_n\to0$, we also have $\varepsilon_n\|u_n-u\|_X \to0$. Taking the limit superior as $n \to \infty$ in the last inequality, we obtain
%Applying Fubini's theorem and H\"older's inequality, we obtain
%\begin{align*}
%|\langle\Phi'_\lambda(u_n), u_n - u \rangle| 
%&\le \|a\|_\infty \int_\Omega |u_n||u_n - u|\,dx + \|K\|_\infty \left( %\int_\Omega |u_n(y)|\,dy \right) \left( \int_\Omega |u_n(x) - u(x)|\,dx \right) \\
%&\quad + \lambda \|Q\|_\infty \left( \int_\Omega |u_n(y)|^q\,dy \right) \left( \int_\Omega |u_n(x)|^{q-1} |u_n(x) - u(x)|\,dx \right) \\
%&\le \|a\|_\infty \|u_n\|_{2} \|u_n - u\|_{2} + \|K\|_\infty \|u_n\|_{1} \|u_n - %u\|_{1} + \lambda \|Q\|_\infty \|u_n\|_{q}^{2q-1} \|u_n - u\|_{q} \to 0
%\end{align*}
\[
\limsup_{n \to \infty} \|u_n\|_X\le \|u\|_X.
\]
On the other hand, the norm is weakly lower semicontinuous, so 
\[
\|u\|_X \le \liminf_{n \to \infty}\|u_n\|_X.
\]
Therefore $\|u_n\|_X \to \|u\|_X$, which implies that $u_n \to u$ strongly in $E$.

\end{proof}

\begin{proposition}\label{CIP}
Assume $(H_1)-(H_3)$. Let $\lambda \ge 0$ and let $\bar{u} \in E$ be a critical point of the functional  $I_{E,\lambda}$ and  
\[
g(\bar u) := a(x)\bar{u}(x)
+ \int_{\Omega} K(x,y)\bar{u}(y)\,dy + \lambda \bar{u}(x)^{q-1}
\int_{\Omega} Q(x,y)|\bar{u}(y)|^q\,dy.
\]
Then the problem 
\begin{equation}\label{dirichet_problem}
\begin{cases}
-\Delta_r \bar{v}=g(\bar u) & \text{in } \Omega\\
\bar{v}=0 & \text{on } \partial\Omega
\end{cases}
\end{equation}
has a unique weak solution $\bar{v}$, and $\bar{v}\in E$. 
\end{proposition}

\begin{proof}
The operator $-\Delta_r:X\to X^*$ is continuous, strictly monotone, and  coercive. Since $\Phi_\lambda\in C^1(X,\mathbb{R})$, Lemma \ref{lem:funcional-c1} shows that the map
\[
\varphi\mapsto\int_\Omega g(\bar u)\varphi\,dx
\]
is a continuous linear functional on $X$.  Hence, by the Minty-Browder theorem \cite[Theorem 5.16]{brezis2011}, there exists a unique $\bar v\in X$ such that
\begin{equation}\label{CIP-weak}
\int_\Omega |\nabla \bar{v}|^{r-2} \nabla \bar{v} \cdot \nabla \varphi \, dx = \int_\Omega g(\bar u)\varphi \, dx \quad \text{for all } \varphi \in X.
\end{equation}

It remains to show that $\bar v \ge 0$ a.e. in $\Omega$. Since $\bar{u} \in E$, $a(x) \ge a_0 > 0$ a.e. in $\Omega$, $K(x,y) \ge 0$ a.e. in $\Omega \times \Omega$, $Q(x,y) \ge q_0 > 0$ a.e. in $\Omega \times \Omega$ and $\lambda \ge 0$, we have 
\[
g(\bar u) \ge 0 \quad \text{a.e. in } \Omega.
\]
Let $\bar{v}^- := \max\{-\bar{v}, 0\} \in X$. Taking $\varphi = \bar{v}^-$ in \eqref{CIP-weak}, we obtain
\[
\int_\Omega |\nabla \bar{v}|^{r-2} \nabla \bar{v} \cdot \nabla \bar{v}^- \, dx = \int_\Omega g(\bar u)\bar{v}^- \, dx.
\]
Since $\nabla \bar{v} = -\nabla \bar{v}^-$ a.e. on $\{\bar{v} < 0\}$ and $\nabla \bar{v}^-=0$ a.e. on $\{\bar{v}> 0\}$, the left-hand side equals
\[
-\int_\Omega |\nabla \bar{v}^-|^r \, dx = -\|\bar{v}^-\|_X^r.
\]
The right-hand side is nonnegative, because $g(\bar u) \ge 0$ and $\bar{v}^- \ge 0$ a.e. in $\Omega$. Hence, 
\[
\|\bar{v}^-\|_X^r \le 0,
\]
so $\bar{v}^- = 0$ a.e. in $\Omega$. That is, $\bar{v} \ge 0$ a.e. in $\Omega$, and therefore $\bar{v} \in E$.

\end{proof}

%\begin{theorem}
%The problem \eqref{eq:problem} admits at least one nontrivial and nonnegative weak solution $u \in X$.
%\end{theorem}

\begin{proof}[Proof of Theorem 1.1]
Fix $\lambda\geq0$. We apply the Mountain Pass Theorem \ref{Mountain_Pass} to the functional
\[
I_{E,\lambda}:=\Psi_E-\Phi_\lambda.
\]
Since $E$ is a closed convex subset of $X$ with $0 \in E$, and $\Psi$ is convex and continuous, $\Psi_E$ is proper, convex and lower
semicontinuous. By Lemma \ref{lem:funcional-c1}, $\Phi_\lambda\in C^1(X,\mathbb{R})$. 

We first verify the Mountain-Pass geometry. Choose any $w\in E\setminus\{0\}$. By $(H_1)$ and $(H_2)$, we have $B(w)>0$. Since $r<2$, 
\[
I_0(tw)
=
\frac{t^r}{r}\|w\|_X^r
-\frac{t^2}{2}B(w)
\to-\infty
\qquad\text{as }t\to+\infty.
\]
Choose $T>0$ sufficiently large and set $e:=Tw\in E$, so that $I_0(e)<0$. Since $C(e)\geq0$ and $\lambda\geq0$, 
\[
I_\lambda(e)
=
I_0(e)-\frac{\lambda}{2q}C(e)
\leq I_0(e)<0.
\]
Thus, $e$ satisfies the negative-energy condition for every $\lambda\geq0$.

By the continuous Sobolev embeddings, there exist constants $S_B,S_C>0$ such that
\[
|B(u)|\leq S_B\|u\|_X^2,
\qquad
0\leq C(u)\leq S_C\|u\|_X^{2q},
\qquad u\in X.
\]
Hence every $u\in E$ satisfies
\[
I_{E,\lambda}(u)
\geq
\frac1r\|u\|_X^r
-\frac{S_B}{2}\|u\|_X^2
-\frac{\lambda S_C}{2q}\|u\|_X^{2q}.
\]
Since $r<2\leq q$, we can choose $\rho \in (0, \|e\|_X)$ sufficiently small that 
\[
\frac{S_B}{2}\rho^{2-r}
+
\frac{\lambda S_C}{2q}\rho^{2q-r}
\leq\frac1{2r}.
\]
Then
\[
I_{E,\lambda}(u)
\geq\frac1{2r}\rho^r=:\alpha_\lambda>0
\qquad
\text{for all } u \in E \text{ with } \|u\|_X=\rho.
\]
For $u \in X \setminus E$ we have $I_{E,\lambda}(u)=+\infty>0$, so $I_{E,\lambda}(u)>0$ for all $u \in X$ with $\|u\|_X=\rho$. Together with $I_{E,\lambda}(0)=0$ and $I_{E,\lambda}(e)<0$, this verifies hypotheses (i)-(iii) of Theorem \ref{Mountain_Pass}.

Now, define
\[
\Gamma_E:=
\left\{
\gamma\in C([0,1],X):
\gamma(0)=0,\ \gamma(1)=e,\
\gamma(t)\in E\ \text{for all }t\in[0,1]
\right\}
\]
and
\[
d_\lambda:=
\inf_{\gamma\in\Gamma_E}
\max_{t\in[0,1]}I_{E,\lambda}(\gamma(t)).
\]
The class $\Gamma_E$ is nonempty, since the path
$\gamma(t)=te$ belongs to it. Every $\gamma\in\Gamma_E$ crosses the sphere with radius $\|u\|_X=\rho$, because $\rho < \|e\|_X$. Therefore,
\[
d_\lambda
\geq
\inf_{\|u\|_X=\rho}I_{E,\lambda}(u)
\geq\alpha_\lambda>0.
\]
Let $\Gamma:= \{\gamma\in C([0,1],X): \gamma(0)=0,\ \gamma(1)=e\}$. Since $\Psi_E=+\infty$ on $X\setminus E$, we have $I_{E,\lambda}=+\infty$ outside $E$. Hence any $\gamma\in\Gamma$ that leaves $E$ gives $\sup_t I_{E,\lambda}(\gamma(t))=+\infty$ and does not affect the infimum. Therefore the minimax value in Theorem~\ref{Mountain_Pass} equals $d_\lambda$, the infimum over $\Gamma_E$.

By Lemma \ref{lem:PS-IE}, the functional $I_{E,\lambda}$ satisfies the Palais-Smale condition at every level. Thus all the hypotheses of Theorem \ref{Mountain_Pass} are fulfilled, and there exists a critical point $\bar u_\lambda\in X$ of $I_{E,\lambda}$ such that
\[
I_{E,\lambda}(\bar u_\lambda)=d_\lambda.
\]
Since $\Psi_E$ is finite only on $E$, we have $\bar u_\lambda\in E$. Also, 
\[
I_\lambda(\bar u_\lambda)
=
I_{E,\lambda}(\bar u_\lambda)
=
d_\lambda
\geq\alpha_\lambda>0
=
I_\lambda(0),
\]
so $\bar u_\lambda\neq0$.

Set 
\[
c_\lambda:=
\inf_{\gamma\in\Gamma}
\max_{t\in[0,1]}I_\lambda(\gamma(t)).
\]
Since $\Gamma_E\subset\Gamma$ and $I_{E,\lambda}=I_\lambda$ on $E$, we obtain
\[
c_\lambda\leq d_\lambda.
\]
Conversely, let $\gamma\in\Gamma$. Since $e \in E$, the path $t \mapsto |\gamma(t)|$ belongs to $\Gamma_E$. Moreover, $I_\lambda(|u|)\leq I_\lambda(u)$ for all $u\in X$, because $C(|u|)=C(u)$ and $B(|u|) \ge B(u)$. Hence 
\[
d_\lambda
\leq
\max_{t\in[0,1]}I_{E,\lambda}(|\gamma(t)|)
\leq
\max_{t\in[0,1]}I_\lambda(\gamma(t)).
\]
Taking the infimum over $\gamma\in\Gamma$ gives $d_\lambda\leq c_\lambda$. Therefore $d_\lambda=c_\lambda$. In particular, 
\[
I_\lambda(\bar u_\lambda)=c_\lambda>0.
\]
It remains to show that $\bar u_\lambda$ is a weak solution of problem \eqref{eq:problem}. Since $\bar u_\lambda\in E$ is a critical point of $I_{E,\lambda}$ and $\lambda\geq0$, Proposition \ref{CIP} gives a function $\bar v_\lambda\in E$ that is a weak solution of the problem
\[
\begin{cases}
-\Delta_r\bar v_\lambda=g(\bar u_\lambda)
&\text{in }\Omega,\\
\bar v_\lambda=0&\text{on }\partial\Omega,
\end{cases}
\]
where
\[
g(\bar u_\lambda)
=a(x)\bar u_\lambda(x)
+\int_\Omega K(x,y)\bar u_\lambda(y)\,dy 
+\lambda|\bar u_\lambda(x)|^{q-2}\bar u_\lambda(x)
\int_\Omega Q(x,y)|\bar u_\lambda(y)|^q\,dy.
\]

By Theorem \ref{adapted_principle}, $\bar u_\lambda$ is a weak solution of problem \eqref{eq:problem}. Since $\bar u_\lambda\in E\setminus\{0\}$, this solution is nonnegative and nontrivial. Therefore, for every $\lambda\geq0$, problem $\eqref{eq:problem}$ has a nontrivial nonnegative weak solution.

\end{proof}

\section{The Competitive Case $\lambda<0$: Variational Thresholds and Multiplicity of Solutions}\label{se.4}

For $\lambda < 0$, problem \eqref{eq:problem} has a logistic-type structure. Here, the two negative variational thresholds allow us to combine global minimization with the Mountain-Pass Theorem, yields multiple solutions.

%Note that, since $a(x) \ge a_0 > 0$ a.e. in $\Omega$ and $K(x,y) \ge 0$ a.e. in $\Omega \times \Omega$, we have 
%\[
%B(u) = \int_\Omega a(x)u^2dx + \iint_{\Omega \times \Omega}K(x,y)u(x)u(y)\,dx\,dy \ge a_0 \int_\Omega u^2 \, dx > 0 \quad \text{for all } u \in E \setminus \{0\}.
%\]
For each $u \in X \setminus \{0\}$, we define the function $\lambda_u: (0,\infty)\to\mathbb{R}$ as the value of $\lambda$ for which $I_{\lambda}(tu)=0$. Solving for $\lambda$ gives
\begin{equation}\label{identity_lambda}
\lambda_u(t)
=
\frac{2q}{C(u)}
\left(
\frac{\|u\|_X^r}{r}t^{r-2q}
-\frac{B(u)}{2}t^{2-2q}
\right).
\end{equation}
Since \(C(u)>0\) for every \(u\in X\setminus\{0\}\), the function \(\lambda_u\) is well defined.

%The derivative of $\lambda_u$ with respect to $t$ is given by
%\begin{equation}\label{derivate_lambda}
%\lambda_u'(t)
%=
%\frac{2q}{C(u)}t^{-2q}
%\left(
%tB(u)(q-1)
%-t^{r-1}\frac{(2q-r)}{r}\|u\|_X^r
%\right).
%\end{equation}
%Setting $\lambda_u'(t_u^*) = 0$, it follows that
For every \(u\in E\setminus\{0\}\), the equation \(\lambda_u'(t)=0\) has a unique positive solution, given by
\begin{equation}\label{critical-points}
t_u^* = \left( \frac{2q-r}{r(q-1)} \frac{\|u\|_X^r}{B(u)} \right)^{\frac{1}{2-r}}.
\end{equation}
It is the global minimum point of $\lambda_u$, and the corresponding minimum value $\lambda_u^* := \min_{t>0}\lambda_u(t) = \lambda_u(t_u^*)$ is
\begin{equation}\label{identity-lambda_u}
\lambda_u^* = -\frac{q(2-r)}{C(u)} \left( \frac{B(u)}{2q-r} \right)^{\frac{2q-r}{2-r}} \left( \frac{r(q-1)}{\|u\|_X^r} \right)^{\frac{2(q-1)}{2-r}}<0.
\end{equation}

Likewise, for each $u \in X \setminus \{0\}$, we define the function $\widehat{\lambda}_u : (0,\infty) \to \mathbb{R}$ as the value of $\lambda$ for which $\langle I_\lambda'(tu), tu \rangle = 0$. Solving for $\lambda$ gives 
\begin{equation}\label{identity_hat-lambda}
\widehat{\lambda}_u(t) := \frac{\|u\|_X^r t^{r-2q} - B(u) t^{2-2q}}{C(u)}.
\end{equation}
An analogous computation shows that \(\widehat{\lambda}_u\) has a unique critical point in \((0,\infty)\) for every \(u\in E\setminus\{0\}\), given by
\[
\hat{t}_u^* = \left( \frac{2q - r}{2(q - 1)} \frac{\|u\|_X^r}{B(u)} \right)^{\frac{1}{2-r}} = \left( \frac{r}{2} \right)^{\frac{1}{2-r}} t_u^*,
\]
which is a multiple of $t_u^*$ that depends only on $r$. The minimum value $\widehat{\lambda}_u^* := \min_{t>0} \widehat{\lambda}_u(t) = \widehat{\lambda}_u(\hat{t}_u^*)$ is 
\begin{equation}\label{identity-hat-lambda_u}
\widehat{\lambda}_u^* = -\frac{2 - r}{C(u)} \left( \frac{B(u)}{2q - r} \right)^{\frac{2q - r}{2 - r}} \left( \frac{2(q - 1)}{\|u\|_X^r} \right)^{\frac{2(q - 1)}{2 - r}} = \frac{1}{q} \left( \frac{2}{r} \right)^{\frac{2q-2}{2-r}} \lambda_u^*<0,
\end{equation}
which is a multiple of $\lambda_u^*$ that depends only on $r$ and $q$.
%\begin{lemma}\label{lem:lambda_u*<0} 
%For all $u \in E$, with $u \neq 0$, such that
%\[
%\lambda_u^*:=\min_{t>0}\lambda_u(t)<0.
%\]
%\end{lemma}

%\begin{proof}
%Let $u \in E$ with $u \neq 0$. By definition \ref{identity_lambda},
%\[
%\lambda_u(t)
%=
%\frac{2q}{C(u)}t^{2-2q}
%\left(
%\frac{\|u\|_X^r}{r}t^{r-2}
%-\frac{B(u)}{2}
%\right).
%\]
%Hence, since $r<2$ and $B(u)>0$ for all $u \in E \setminus \{0\}$, we obtain that
%\[
%\frac{\|u\|_X^r}{r}t^{r-2}
%-\frac{B(u)}{2} 
%\to
%-\frac{B(u)}{2}<0, \qquad \text{as } t \to +\infty.
%\] 
%Therefore, there exists $T>0$ such that
%\[
%\frac{\|u\|_X^r}{r}t^{r-2}
%-\frac{B(u)}{2}<0
%\]
%for all $t>T$. Since
%\[
%\frac{2q}{C(u)}t^{2-2q}>0,
%\]
%it follows that
%\[
%\lambda_u(t)<0
%\qquad\text{for all } t>T.
%\]
%Therefore,
%\[
%\min_{t>0}\lambda_u(t)<0, \quad \text{ for all } u \in E.
%\]
%\end{proof}

\begin{lemma}\label{lem:lambda*} 
Assume $(H_1)-(H_3)$. Then the infimum 
\[
\lambda^* := \inf_{u\in X\setminus\{0\}} \inf_{t>0}\lambda_u(t)
\]
is attained at some $u^* \in E \setminus \{0\}$ and satisfies $\lambda^* > -\infty$.
\end{lemma}

\begin{proof}
First, for any $u \in X \setminus \{0\}$, taking $|u| \in E \setminus \{0\}$, we have $\||u|\|_X=\|u\|_X$, $C(|u|) = C(u)$ and $B(|u|) \ge B(u)$. By the identity of $\lambda_u(t)$ in \eqref{identity_lambda}, it follows that $\min_{t>0}\lambda_{|u|}(t) \le \inf_{t>0}\lambda_u(t)$. Consequently, 
\[
\lambda^* = \inf_{u\in X\setminus\{0\}} \inf_{t>0}\lambda_u(t) =\inf_{u\in E\setminus\{0\}} \min_{t>0}\lambda_u(t). 
\]
Now, for all $u\in E\setminus\{0\}$ and $s>0$, we obtain
\[
\lambda_{su}^*
=
\min_{t>0}\lambda_{su}(t)
=
\min_{t>0}\lambda_u(st)
=
\min_{k>0}\lambda_u(k)
=
\lambda_u^*,
\]
which shows that $\lambda_u^*$ is scale invariant. Taking $s=C(u)^{-1/(2q)}$, we have $s u \in E$ and $C(su)=s^{2q}C(u)=1$. By the scale invariance of $\lambda_u^*$, it follows that
\[
\lambda^*
=
\inf_{u\in M}\lambda_u^*,
\]
where $M:=\{u\in E : C(u)=1\}$.

We divide the remaining proof into two claims. \\
\textbf{Claim 1:} $\lambda^* > -\infty$. Indeed, for any $u \in M$, since $Q(x,y) \ge q_0 > 0$ a.e. in $\Omega \times \Omega$, we have
\[
1 = C(u) \ge q_0 \iint_{\Omega\times\Omega} |u(x)|^q |u(y)|^q \, dx \, dy = q_0 \|u\|_{q}^{2q},
\]
which implies $\|u\|_{q} \le q_0^{-1/(2q)}$. On the other hand, since $q \ge 2$, by continuous embeddings in Lebesgue spaces, there exist constants $C_1, C_2 > 0$ such that
\begin{equation}\label{4.2-1}
\begin{aligned}
|B(u)| 
&\le \int_\Omega |a(x)| u^2 \, dx + \iint_{\Omega\times\Omega} |K(x,y)| |u(x)| |u(y)| \, dx \, dy \\
&\le \|a\|_\infty \|u\|_{2}^2 + \|K\|_\infty \|u\|_{1}^2 \\
&\le \left( C_1 \|a\|_\infty + C_2 \|K\|_\infty \right) \|u\|_{q}^2 \\
&\le \left( C_1 \|a\|_\infty + C_2 \|K\|_\infty \right) q_0^{-1/q} =: \alpha.
\end{aligned}
\end{equation}
In addition, using the Sobolev embedding $\|u\|_q\le C_S\|u\|_X$, we obtain
\begin{equation}\label{4.2-2}
1 = C(u) \le \|Q\|_\infty \|u\|_{q}^{2q} \le C_S^{2q}\|Q\|_{\infty} \|u\|_X^{2q} \implies \|u\|_X \ge \beta > 0 \quad \forall u \in M,
\end{equation}
where $\beta := C_S^{-1}\|Q\|_{\infty}^{-\frac{1}{2q}} > 0$. Substituting \eqref{4.2-1} and \eqref{4.2-2} into the identity \eqref{identity-lambda_u}, we obtain
\[
-\lambda_u^* = q(2-r) \left( \frac{B(u)}{2q-r} \right)^{\frac{2q-r}{2-r}} \left( \frac{r(q-1)}{\|u\|_X^r} \right)^{\frac{2(q-1)}{2-r}} \le q(2-r) \left( \frac{\alpha}{2q-r} \right)^{\frac{2q-r}{2-r}} \left( \frac{r(q-1)}{\beta^r} \right)^{\frac{2(q-1)}{2-r}} =: C_0,
\]
for all $u \in M$, which implies that $\lambda_u^* \ge -C_0 > -\infty$. Taking the infimum over $M$, we conclude that $\lambda^* \ge -C_0 > -\infty$. \\
\textbf{Claim 2:} The infimum defining $\lambda^*$ is attained at some $u^* \in M$. Indeed, let $(u_n) \subset M$ be a minimizing sequence, that is, $\lambda_{u_n}^* \to \lambda^*$. For each $n$, let $t_n > 0$ be the minimizer of $t \mapsto \lambda_{u_n}(t)$, so that $\lambda_{u_n}(t_n) = \lambda_{u_n}^*$.

Note that $\lambda_u^* < 0$ for all $u \in E \setminus \{0\}$, which ensures that $\lambda^* \le \lambda_u^* < 0$. Since $\lambda_{u_n}^* \to \lambda^*$, there exist $n_0 \in \mathbb{N}$ and $c_1 > 0$ such that
\begin{equation}\label{4.2-3}
c_1 \le -\lambda_{u_n}^* = q(2-r) \left( \frac{B(u_n)}{2q-r} \right)^{\frac{2q-r}{2-r}} \left( \frac{r(q-1)}{\|u_n\|_X^r} \right)^{\frac{2(q-1)}{2-r}} \qquad \forall n \ge n_0.
\end{equation}
Isolating $\|u_n\|_X$ in \eqref{4.2-3} and using $|B(u_n)| \le \alpha$ established in \eqref{4.2-1}, we obtain
\begin{equation}\label{4.2-4}
\|u_n\|_X^{\frac{2r(q-1)}{2-r}} \le \frac{q(2-r)}{c_1} \left( \frac{\alpha}{2q-r} \right)^{\frac{2q-r}{2-r}} [r(q-1)]^{\frac{2(q-1)}{2-r}},
\end{equation}
which proves that $(u_n)$ is bounded in $X$. Since $\|u_n\|_X \ge \beta>0$ by \eqref{4.2-2}, inequality \eqref{4.2-3} gives 
\[
c_1 \le q(2-r)\left(\frac{B(u_n)}{2q-r}\right)^{\frac{2q-r}{2-r}}\left(\frac{r(q-1)}{\beta^r}\right)^{\frac{2(q-1)}{2-r}} \quad \forall n \ge n_0, 
\]
so $B(u_n)$ is bounded away from zero, that is, there exists $c_2 > 0$ such that
%Consequently, since $-\lambda_{u_n}^* \ge c_1 > 0$ in \eqref{4.2-3} and $\|u_n\|_X$ is bounded from above in \eqref{4.2-4}, it follows that $B(u_n)$ is bounded away from zero, i.e., there exists $c_2 > 0$ such that
\begin{equation}\label{4.2-5}
B(u_n) \ge c_2 > 0 \qquad \forall n \ge n_0.
\end{equation}
Since $\|u_n\|_X$ and $B(u_n)$ are both bounded and bounded away from zero, the sequence $(t_n)$ given in \eqref{critical-points} by 
\[
t_n = \left( \frac{2q-r}{r(q-1)} \frac{\|u_n\|_X^r}{B(u_n)} \right)^{\frac{1}{2-r}}
\]
satisfies $0 < \delta_1 \le t_n \le \delta_2 < +\infty$ for all $n \ge n_0$.

From the reflexivity of $X$ and the boundedness of the sequence $(u_n)$ in $X$, there exists a subsequence (still denoted by $u_n$) and $u^* \in X$ such that
\begin{align*}
u_n &\rightharpoonup u^* \,\,\text{weakly in } X. 
%u_n &\to u^* \quad \text{strongly in } L^s(\Omega) \text{ for all } 1 \le s < r^*.
\end{align*}
Since $E$ is closed and convex in $X$, it is weakly closed; thus $(u_n) \subset E$ implies $u^* \in E$. By Lemma \ref{lem:convergences}, it follows that $C(u_n) \to C(u^*)$ and $B(u_n) \to B(u^*)$. Since $C(u_n) = 1$, we have $C(u^*) = 1$, which together with $u^* \in E$ shows that $u^* \in M$. %$u_n \rightharpoonup u^*$ weakly in $X$. Since $E$ is closed and convex in $X$, it is weakly closed; thus, since $(u_n) \subset E$, we have $u^* \in E$. By Lemma \ref{lem:convergences}, it follows that $C(u_n) \to C(u^*)$ and $B(u_n) \to B(u^*)$. Since $C(u_n) = 1$, we have $C(u^*) = 1$, which together with $u^* \in E$ shows that $u^* \in M$. 

Passing to a subsequence if necessary, we may assume that $t_n \to t_* \in [\delta_1, \delta_2]$. By the lower semicontinuity of the norm, $\|u^*\|_X \le \liminf_{n\to\infty} \|u_n\|_X$. Taking the limit inferior in \eqref{identity_lambda}, evaluated at $u_n$ and $t_n$, and using $C(u_n)=C(u^*)=1$, $B(u_n) \to B(u^*)$ and $t_n \to t_*$, we obtain
\[
\lambda^*= \liminf_{n\to\infty} \lambda_{u_n}^*=\liminf_{n\to\infty} \lambda_{u_n}(t_n) \ge \lambda_{u^*}(t_*).
\]
On the other hand, by definition of $\lambda^*$ and since $u^* \in M$, we have $\lambda^* \le \lambda_{u^*}^* \le \lambda_{u^*}(t_*)$. It follows that
\[
\lambda_{u^*}(t_*) = \lambda_{u^*}^* = \lambda^*,
\]
which proves that $\lambda^*$ is attained by $u^* \in M$. This completes the proof. 

\end{proof}

\begin{lemma}\label{lem:hat-lambda*}
Assume $(H_1)-(H_3)$. Then the infimum
\[
\widehat{\lambda}^* := \inf_{u\in X\setminus\{0\}} \inf_{t>0}\widehat{\lambda}_u(t)
\]
is attained and satisfies $\widehat{\lambda}^* > -\infty$. Moreover, $\widehat{\lambda}^*$ is attained at every $u \in E \setminus \{0\}$ at which $\lambda^*$ is attained.
\end{lemma}

\begin{proof}
Recall that for every $u \in E \setminus \{0\}$, the unique minimum $\widehat{\lambda}_u^* := \min_{t>0} \widehat{\lambda}_u(t)$ satisfies
\[
\widehat{\lambda}_u^* = \kappa \lambda_u^*, \quad \text{where } \kappa := \frac{1}{q} \left( \frac{2}{r} \right)^{\frac{2(q-1)}{2-r}} > 0.
\]
Now let $u \in X \setminus \{0\}$. Since $\||u|\|_X=\|u\|_X$, $C(|u|)=C(u)$, and $B(|u|) \ge B(u)$, we obtain
\[
\inf_{t>0} \widehat{\lambda}_u(t) \ge \widehat{\lambda}_{|u|}^*,
\]
where $|u| \in E \setminus \{0\}$. Together with $E \setminus \{0\} \subset X \setminus \{0\}$, this yields
\[
\widehat{\lambda}^* = \inf_{u \in X \setminus \{0\}} \inf_{t>0} \widehat{\lambda}_u(t)= \inf_{u \in E \setminus \{0\}} \widehat{\lambda}_u^* = \kappa \inf_{u \in E \setminus \{0\}} \lambda_u^* = \kappa \lambda^*.
\]
By Lemma \ref{lem:lambda*}, $\lambda^*>-\infty$ and $\lambda^*$ is attained at some $u^* \in E \setminus \{0\}$, that is, $\lambda_{u^*}^* = \lambda^*$. Hence, $\widehat{\lambda}^* = \kappa \lambda^* > -\infty$ and $\widehat{\lambda}_{u^*}^* = \kappa \lambda_{u^*}^* = \kappa \lambda^* = \widehat{\lambda}^*$, so $\widehat{\lambda}^*$ is attained at $u^*$.

\end{proof}
\begin{remark}
We have $\widehat{\lambda}^*<\lambda^*$.
Since $\widehat{\lambda}^*=\kappa\lambda^*$ and
$\lambda^*<0$, it suffices to show that $\kappa>1$.
Setting $s=2/r>1$ and using
$\log s>1-1/s$ and $\log q\leq q-1$, we obtain
\[
\log\kappa
=(q-1)\frac{s\log s}{s-1}-\log q
>q-1-\log q
\geq0.
\]
Hence $\kappa>1$, as required.
\end{remark}

\subsection*{Coexistence of Positive- and Negative-Energy Solutions for $\lambda \in (\lambda^*, 0)$}

To obtain two distinct critical points with opposite energy signs, it is necessary that $I_\lambda$ take negative values, which does not occur for every $\lambda<0$. For each $\lambda \in (\lambda^*, 0)$, global minimization yields a critical point with negative energy, whereas the Mountain-Pass Theorem provides a critical point at positive energy. 
%The variational geometry required for the coexistence of two distinct critical points with opposite energy signs is not automatically available for every $\lambda<0$. Indeed, for certain negative values of $\lambda$, the energy functional may remain strictly positive along all fibers, so no negative energy levels exist. So we look for an interval of negative values of $\lambda$ on which the functional attains negative values. On such an interval the geometry required by the Mountain-Pass Theorem is available, and combining it with global minimization gives two critical points whose energy levels have opposite signs, hence two distinct solutions to problem \eqref{eq:problem}. Below we establish the conditions under which both methods apply.

Since $\lambda^*$ is attained, there exists a direction along which $I_\lambda$ is nonpositive for every $\lambda \in [\lambda^*,0)$. Indeed, by Lemma \ref{lem:lambda*}, there exist $u_* \in X \setminus \{0\}$ and $t_* > 0$ such that $\lambda^* = \lambda_{u_*}(t_*)$. By the definition of $\lambda_u(t)$ in \eqref{identity_lambda}, we have $I_{\lambda^*}(t_* u_*) = 0$. Hence,
\[
I_\lambda(t_* u_*) = I_{\lambda^*}(t_* u_*) - \frac{\lambda - \lambda^*}{2q} C(t_*u_*) \le 0.
\]

\begin{lemma}
\label{lem:coercive}
Assume $(H_1)-(H_3)$. Then, for every $\lambda<0$, the functional $I_\lambda:X\to\mathbb{R}$ is coercive in $X$.
\end{lemma}

\begin{proof}
Since $q\ge 2$, it follows from the continuous embeddings that there exists $\overline{C}>0$ such that
\[
\left|B(u)\right|
\leq \|a\|_\infty\|u\|_2^2 + \|K\|_{\infty}\|u\|_1^2
\leq \overline{C}\|u\|_q^2.
\]
Since $Q(x,y)\geq q_0$ a.e. in $\Omega\times\Omega$, we also have 
\[
C(u)
\geq q_0
\iint_{\Omega\times\Omega}
|u(x)|^q|u(y)|^q\,dx\,dy =q_0\|u\|_q^{2q}.
\]
Combining these estimates and using $\lambda<0$, we obtain constants $C_1,C_2>0$ such that
\[
-\frac{\lambda}{2q}C(u) 
-\frac12B(u)
\geq
C_1\|u\|_q^{2q}-C_2\|u\|_q^2.
\]
The function
\[
g(s)=C_1s^{2q}-C_2s^2,\qquad s\geq0,
\]
is continuous. Since $q\geq2$, we have $\lim_{s\to+\infty}g(s)=+\infty$, while $g(s)<0$ for all sufficiently small $s>0$. Thus, $g$ attains a finite, strictly negative global minimum on $[0,+\infty)$. In particular, there exists $M>0$ such that
\[
g(s)\geq-M,\qquad s\geq0.
\]
It follows that
\[
C_1\|u\|_q^{2q}-C_2\|u\|_q^2\geq -M, 
\] 
and hence, 
\[
I_\lambda(u)
\geq
\frac{1}{r}\|u\|_X^r-M.
\]
Since $r>0$, we obtain 
\[
I_\lambda(u)\to+\infty
\qquad\text{as }\|u\|_X\to\infty.
\]
Therefore, $I_\lambda$ is coercive in $X$.

\end{proof}

\begin{proposition}[Existence of a nontrivial global minimizer]\label{prop:minimizer}
Assume $(H_1)-(H_3)$. For each $\lambda \in (\lambda^*, 0)$, the functional $I_{\lambda}$ admits a nontrivial global minimizer $v_\lambda \in X \setminus \{0\}$ satisfying
\[
I_{\lambda}(v_\lambda) = \inf_{v\in X} I_{\lambda}(v) < 0.
\]
\end{proposition}

\begin{proof}
Fix $\lambda \in (\lambda^*, 0)$ and let $(v_n) \subset X$ be a minimizing sequence for $I_\lambda$, so that
$$
I_\lambda(v_n) \to \inf_{v\in X} I_\lambda(v).
$$
By Lemma \ref{lem:coercive}, $I_\lambda$ is coercive in $X$, which implies that the sequence $(v_n)$ is bounded in $X$. Since $X $ is a reflexive Banach space, after passing to a subsequence, there exists $v_\lambda \in X$ such that
$$
v_n \rightharpoonup v_\lambda \quad \text{weakly in } X.
$$
By the weak lower semicontinuity of the norm and Lemma \ref{lem:convergences}, the functional $I_\lambda$ is weakly lower semicontinuous in $X$. Therefore,
$$
I_\lambda(v_\lambda) \le \liminf_{n\to\infty} I_\lambda(v_n) = \inf_{v\in X} I_\lambda(v),
$$
so $v_\lambda \in X$ is a global minimizer of $I_\lambda$. 

It remains to show that the minimum is negative. Since $\lambda \in (\lambda^*, 0)$, there exist $t_*>0$ and $v_* \in X \setminus \{0\}$ such that $I_\lambda(t_* v_*) < 0$. Consequently,
$$
I_\lambda(v_\lambda) = \inf_{v\in X} I_\lambda(v) \le I_\lambda(t_*v_*) < 0.
$$
Since $I_\lambda(0)=0$, the minimizer $v_\lambda$ is nontrivial. 

\end{proof}

\begin{lemma}\label{lem:PS_c}
Assume $(H_1)-(H_3)$. Then, for every $\lambda<0$ and every $c \in \R$, the functional $I_{\lambda}$ satisfies the Palais-Smale condition at level $c$.
\end{lemma}

\begin{proof} 
Let $(u_n) \subset X$ be a Palais-Smale sequence at level $c$. By Lemma \ref{lem:coercive}, $I_\lambda$ is coercive, so the sequence $(u_n)$ is bounded in $X$. By reflexivity, up to a subsequence, there exists $u \in X$ such that
\[
u_n \rightharpoonup u \text{ weakly in } X.
\]
Because $I_\lambda'(u_n) \to 0$ in $X^*$, we have $\langle I_\lambda'(u_n), u_n - u \rangle \to 0$. By Lemma \ref{lem:convergences}, we obtain
\[
\int_\Omega |\nabla u_n|^{r-2}\nabla u_n \nabla(u_n - u)\,dx \to 0.
\]
Since the $r$-Laplacian operator satisfies the $(S_+)$ condition in $X$, we conclude that $u_n \to u$ strongly in $X$. 

\end{proof}

\begin{lemma}\label{MountainPass_geometry}
Assume $(H_1)-(H_3)$. Then, for each $\lambda \in [\lambda^*, 0)$, the energy functional $I_{\lambda}$ satisfies the Mountain-Pass geometry in $X$.
\end{lemma}

\begin{proof} Clearly, $I_\lambda (0)=0$. Since $C(u) \ge 0$ for all $u \in X$ and $\lambda<0$, the estimate
\[
\left|B(u)\right| \le S_B\|u\|_X^2
\]
for some constant $S_B>0$ gives 
\[
I_{\lambda}(u)=
\frac1r\|u\|_X^r-\frac{1}{2}B(u)
-\frac{\lambda}{2q}C(u) 
\geq
\frac1r\|u\|_X^r
-S_B\|u\|_X^2.
\]
Choosing $\rho>0$ sufficiently small such that 
\[
S_B\rho^{2-r}\leq\frac{1}{2r},
\]
we obtain
\[
I_\lambda(u)\geq \frac{1}{2r}\|u\|_X^r>0, \quad\text{for all }u\in X
\text{ with }0<\|u\|_X \le \rho.
\]
Setting $\alpha:=\rho^r/(2r)$ yields
\[
I_\lambda(u)\geq\alpha>0
\qquad\text{whenever }\|u\|_X=\rho.
\]

For the second geometric condition, by Lemma \ref{lem:lambda*} and since $\lambda\in[\lambda^*,0)$, there exist $u_*\in E \setminus \{0\}$ and $t_*>0$ such that $e:=t_*u_*$ satisfies $I_\lambda(e) \le 0$. Since $I_\lambda(u)>0$ for all $u\in X$ with $0<\|u\|_X \le \rho$, we have $\|e\|_X>\rho$. Therefore, $I_{\lambda}$ satisfies the geometry of the Mountain-Pass Theorem.

\end{proof}

\begin{proposition}\label{MountainPass-Classic}
Assume $(H_1)-(H_3)$ and let \(\lambda\in[\lambda^*,0)\). Then there exists $e \in X \setminus \{0\}$ with $I_\lambda(e)\le0$ such that the functional \(I_\lambda\) admits a critical point \(u_\lambda\in X\) at the Mountain-Pass level \(c_\lambda>0\), where
\[
I_\lambda(u_\lambda)=c_\lambda
\qquad\text{and}\qquad
I_\lambda'(u_\lambda)=0,
\]
with
\[
c_\lambda=\inf_{\gamma\in\Gamma}\max_{t\in[0,1]}I_\lambda(\gamma(t)),
\]
and
\[
\Gamma=\{\gamma\in C([0,1],X):\gamma(0)=0,\ \gamma(1)=e\}.
\]
\end{proposition}

\begin{proof}
Fix \(\lambda\in[\lambda^*,0)\). By Lemma \ref{MountainPass_geometry}, there exist \(\rho,\alpha>0\) and \(e \in X \setminus \{0\}\) such that
\[
I_\lambda(u)\geq\alpha>0
\quad\text{if }\|u\|_X=\rho,
\qquad
I_\lambda(e)\le 0,\quad \|e\|_X>\rho.
\]
Thus, every path \(\gamma\in\Gamma\) crosses the sphere \(\|u\|_X=\rho\), and 
\[
c_\lambda=\inf_{\gamma\in\Gamma}\max_{t\in[0,1]}I_\lambda(\gamma(t))
\geq\alpha>0.
\]
Since \(I_\lambda\in C^1(X,\mathbb{R})\) and satisfies the Palais–Smale condition at level \(c_\lambda\) (Lemma \ref{lem:PS_c}), the Mountain Pass Theorem guarantees the existence of a critical point \(u_\lambda\in X\) such that
\[
I_\lambda(u_\lambda)=c_\lambda>0
\quad\text{and}\quad
I_\lambda'(u_\lambda)=0.
\]

\end{proof}

\begin{proof}[Proof of Theorem 1.2 \rm (i)]
Fix $\lambda \in (\lambda^*, 0)$. By Proposition \ref{prop:minimizer}, $I_\lambda$ has a global minimizer $v_\lambda \in X \setminus \{0\}$ with negative energy. Since $I_\lambda \in C^1(X,\R)$, we conclude that $I_\lambda'(v_\lambda)=0$. By Proposition \ref{MountainPass-Classic}, we have a Mountain-Pass critical point $u_\lambda \in X \setminus \{0\}$ with positive energy. Since their energy levels are distinct, we obtain two distinct, nontrivial solutions to problem \eqref{eq:problem}. 

\end{proof}

\subsection*{At $\lambda = \lambda^*$: Two Distinct Solutions with Zero and Positive Energy}

At the energy threshold $\lambda = \lambda^*$, the energy functional satisfies $I_{\lambda^*}(u) \ge 0$ for all $u \in X$, and hence attains its minimum value of zero. The existence of a nontrivial global minimizer $v_{\lambda^*} \in X \setminus \{0\}$ allows us to obtain a second critical point via the mountain-pass theorem.

%The existence of a nontrivial global minimizer $v_{\lambda^*} \in X \setminus \{0\}$ allows us to obtain multiple nontrivial solutions to \eqref{eq:problem} for $\lambda = \lambda^*$.

%\begin{proposition}\label{lambda=lamda^*}
%For $\lambda = \lambda^*$, the energy functional $I_{\lambda^*}$ admits at least two distinct nontrivial critical points in $X$.
%\end{proposition}

\begin{proof}[Proof of Theorem 1.2 \rm (ii)]
By the definition of $\lambda_v(t)$ in \eqref{identity_lambda}, we have
\begin{equation}\label{eq:fiber-identity-lambda-star}
I_{\lambda^*}(tv)
=
\frac{C(v)}{2q}t^{2q}
\bigl(\lambda_v(t)-\lambda^*\bigr)
\end{equation}
for every $v\in X\setminus\{0\}$ and $t>0$.
Since $\lambda_v(t)\geq\lambda^*$, it follows that
$I_{\lambda^*}(u)\geq0$ for all $u\in X$. Hence
\[
\inf_{u\in X}I_{\lambda^*}(u)=0.
\]
By Lemma \ref{lem:lambda*}, there exist $v_*\in X\setminus\{0\}$ and $t_*>0$ such that $\lambda_{v_*}(t_*)=\lambda^*$. Setting
$v_{\lambda^*}:=t_*v_*$, identity \eqref{eq:fiber-identity-lambda-star} yields
\[
I_{\lambda^*}(v_{\lambda^*})=0
=\inf_{u\in X}I_{\lambda^*}(u).
\]
Thus, $v_{\lambda^*}$ is a nontrivial global minimizer. Since $I_{\lambda^*}\in C^1(X,\mathbb R)$, we also have
\[
I_{\lambda^*}'(v_{\lambda^*})=0.
\]

To obtain a second nontrivial critical point, we apply Lemma \ref{MountainPass_geometry} with \(\lambda=\lambda^*\), which gives
\[
I_{\lambda^*}(u)\geq \frac{1}{2r}\|u\|_X^r>0 \qquad\text{for all }0<\|u\|_X\le\rho.
\]
In particular, 
\[
I_{\lambda^*}(u)\geq\alpha>0
\qquad\text{for all }\|u\|_X=\rho,
\]
where $\alpha:=\rho^r/(2r)$. Taking $e:=v_{\lambda^*}$, define
\[
\Gamma:=
\{\gamma\in C([0,1],X): \gamma(0)=0,\ \gamma(1)=e\}
\]
and
\[
c_{\lambda^*}:=
\inf_{\gamma\in\Gamma}
\max_{t\in[0,1]}I_{\lambda^*}(\gamma(t)).
\]
Since $\|e\|_X > \rho$, every path in $\Gamma$ crosses the sphere
$\|u\|_X=\rho$. Consequently
\[
c_{\lambda^*}\geq
\inf_{\|u\|_X=\rho}I_{\lambda^*}(u)
\geq\alpha>0.
\]

Finally, since $\lambda^*<0$, the Palais-Smale condition holds for $\lambda=\lambda^*$ by Lemma \ref{lem:PS_c}. Therefore, by the Mountain-Pass Theorem, there exists $u_{\lambda^*}\in X$ such that
\[
I_{\lambda^*}'(u_{\lambda^*})=0,
\qquad
I_{\lambda^*}(u_{\lambda^*})=c_{\lambda^*}>0.
\]
In particular, $u_{\lambda^*}\neq0$ and
\[
I_{\lambda^*}(u_{\lambda^*})
>
I_{\lambda^*}(v_{\lambda^*}),
\]
so that $u_{\lambda^*}\neq v_{\lambda^*}$. Thus, $I_{\lambda^*}$ admits at least two distinct nontrivial critical points in $X$.

\end{proof}

\subsection*{Two Distinct Positive-Energy Solutions Below the Energy Threshold $\lambda^*$}\label{sec:two-solutions-below-threshold}

Although $\lambda^*$ separates the parameter range in which negative energy values occur from the range in which they do not, it is not a threshold below which \eqref{eq:problem} has no nontrivial solutions. We show that two nontrivial solutions, both with positive energy, appear immediately below $\lambda^*$.

Recall that $I_{\lambda^*}$ is coercive and admits a nonzero global minimizer $v_{\lambda^*}$ satisfying
\[
I_{\lambda^*}(v_{\lambda^*})=0,
\qquad
I_{\lambda^*}(u)\geq0
\quad\text{for every }u\in X.
\]
For $\lambda<\lambda^*$, the identity
\begin{equation}\label{eq:below-threshold-comparison}
I_\lambda(u)
=
I_{\lambda^*}(u)
+\frac{\lambda^*-\lambda}{2q}C(u)
\end{equation}
shows that $I_\lambda(u)>0$ for every $u\neq0$, since $C(u)>0$. Nevertheless, the energy of $v_{\lambda^*}$ remains close to zero when $\lambda$ is close to $\lambda^*$. We use this to construct a nonzero local minimum and then a second critical point at a higher mountain-pass level. 

%\begin{proposition}\label{prop:two-solutions-below-threshold}
%There exists $\varepsilon>0$ such that, for every $\lambda\in(\lambda^*-\varepsilon,\lambda^*)$, problem \eqref{eq:problem} admits two distinct
%nontrivial weak solutions $v_\lambda,u_\lambda\in X$ with
%\[
%0<I_\lambda(v_\lambda)<I_\lambda(u_\lambda).
%\]
%The solution $v_\lambda$ is a local minimizer of
%$I_\lambda$. Also, $\varepsilon$ can be chosen so that $\lambda^*-\varepsilon>\widehat\lambda^*$.
%\end{proposition}

\begin{proof}[Proof of Theorem 1.2 \rm (iii)] Recall that there exists a constant $S_B>0$ such that
\[
|B(u)|\leq S_B\|u\|_X^2, \quad \text{ for all } u \in X.
\]
Since $r<2$, we can choose $\rho>0$ such that
\[
\rho<\|v_{\lambda^*}\|_X,
\qquad
\frac{S_B}{2}\rho^2\leq\frac{\rho^r}{2r}.
\]
Set $\alpha:=\rho^{r}/(2r)>0$. Since $\lambda<\lambda^*<0$ and $C(u) \ge 0$ for all $u \in X$, we have 
\[
I_\lambda(u)
\geq
\frac1r\|u\|_X^r-\frac{S_B}{2}\|u\|_X^2.
\]
In particular,
\begin{equation}\label{eq:below-threshold-inner-barrier}
I_\lambda(u)\geq\alpha
\qquad\text{if }\|u\|_X=\rho.
\end{equation}
Next, since $I_{\lambda^*}$ is coercive, we can choose $R>\|v_{\lambda^*}\|_X$ sufficiently large such that
\[
I_{\lambda^*}(u)\geq\alpha
\qquad\text{whenever }\|u\|_X\geq R.
\]
Identity \eqref{eq:below-threshold-comparison}
then implies
\begin{equation}\label{eq:below-threshold-outer-barrier}
I_\lambda(u)\geq\alpha
\qquad
\text{if }\|u\|_X=R,\quad \lambda<\lambda^*.
\end{equation}
Choose $\varepsilon>0$ satisfying
\[
\varepsilon<\lambda^*-\widehat\lambda^*,
\qquad
\frac{\varepsilon}{2q}C(v_{\lambda^*})<\alpha,
\]
and fix
$\lambda\in(\lambda^*-\varepsilon,\lambda^*)$.
Since $I_{\lambda^*}(v_{\lambda^*})=0$,
\[
0<I_\lambda(v_{\lambda^*})
=
\frac{\lambda^*-\lambda}{2q}C(v_{\lambda^*})<\alpha.
\]
Thus, $I_\lambda(v_{\lambda^*})$ is strictly below the values of $I_\lambda$ on both boundary spheres.

We now minimize over the closed annular set
\[
S_{\rho,R}
:=
\{u\in X:\rho\leq\|u\|_X\leq R\},
\qquad
m_\lambda:=\inf_{u\in S_{\rho,R}}I_\lambda(u).
\]
Since $v_{\lambda^*}\in S_{\rho,R}$,
\begin{equation}\label{eq:below-threshold-local-level}
m_\lambda\leq I_\lambda(v_{\lambda^*})<\alpha.
\end{equation}

Let $(v_n)\subset S_{\rho,R}$ be a minimizing
sequence. Since $(v_n)$ is bounded in $X$ and $X$ is reflexive,
there exist a subsequence, still denoted by $(v_n)$, and
$v_\lambda\in X$ such that
\[
v_n\rightharpoonup v_\lambda
\quad\text{in }X,
\qquad
v_n\to v_\lambda
\quad\text{in }L^s(\Omega) \text{ for all }\,  1\le s <r^*.
\]
Consequently,
\begin{equation}\label{convergence}
B(v_n)\to B(v_\lambda),
\qquad
C(v_n)\to C(v_\lambda).    
\end{equation}

Although $S_{\rho,R}$ need not be weakly closed, its minimizing sequence cannot converge weakly into the inner ball. Indeed, suppose that $\|v_\lambda\|_X\leq\rho$. Since $\|v_n\|_X\geq\rho$ for every $n$, the convergences in \eqref{convergence} give
\[
\begin{aligned}
m_\lambda
&\geq
\frac{\rho^r}{r}
-\frac12B(v_\lambda)
-\frac{\lambda}{2q}C(v_\lambda)\\
&\geq
\frac{\rho^r}{r}-\frac{S_B}{2}\rho^2
\geq\alpha,
\end{aligned}
\]
contradicting \eqref{eq:below-threshold-local-level}. Hence $\|v_\lambda\|_X>\rho$. 

The weak lower semicontinuity of the norm, together with the convergence of $B$ and $C$, yields
\[
\|v_\lambda\|_X\leq R,
\qquad
I_\lambda(v_\lambda)\leq m_\lambda.
\]
Therefore $v_\lambda\in S_{\rho,R}$ and $I_\lambda(v_\lambda)=m_\lambda$.
By \eqref{eq:below-threshold-outer-barrier}, the strict inequality $m_\lambda<\alpha$ excludes
$\|v_\lambda\|_X=R$. Thus
\[
\rho<\|v_\lambda\|_X<R.
\]
The minimizer lies in the interior of $S_{\rho,R}$ and is therefore a local minimizer in the full space $X$. In particular,
\[
I_\lambda'(v_\lambda)=0,
\qquad
0<I_\lambda(v_\lambda)<\alpha.
\]

By \eqref{eq:below-threshold-comparison}, $I_\lambda(u)\geq I_{\lambda^*}(u)$ for all $u \in X$, so $I_\lambda$ is coercive. Hence every Palais-Smale sequence $(z_n)$ is bounded in $X$. Since $X$ is reflexive, up to a subsequence, there exists $z \in X$ such that
\[
z_n\rightharpoonup z \quad \text{ weakly in } X.
\]
Testing $I_\lambda'(z_n)$ against $z_n-z$ and using Lemma \ref{lem:convergences}, we obtain
\[
\int_\Omega |\nabla z_n|^{r-2}\nabla z_n \nabla (z_n-z)\, dx \to 0.
\]
By the $(S_+)$ property of the $r$-Laplacian operator, we conclude that $z_n \to z$ strongly in $X$. 

We finally apply the mountain-pass theorem in its two-endpoint minimax form. Define
\[
\Gamma_\lambda
:=
\{\gamma\in C([0,1],X):
\gamma(0)=0,\ \gamma(1)=v_\lambda\},
\]
and
\[
c_\lambda
:=
\inf_{\gamma\in\Gamma_\lambda}
\max_{t\in[0,1]}I_\lambda(\gamma(t)).
\]
The class is nonempty, since the path $t\mapsto tv_\lambda$ belongs to $\Gamma_\lambda$. In addition, every path $\gamma\in\Gamma_\lambda$ crosses the sphere $\|u\|_X=\rho$, since $\|v_\lambda\|_X>\rho$. Hence,
\[
c_\lambda\geq\alpha>I_\lambda(v_\lambda)>I_\lambda(0)=0.
\]
The mountain-pass theorem and the Palais-Smale condition yield a critical point $u_\lambda$ at level $c_\lambda$. Thus $u_\lambda$ and $v_\lambda$ are nontrivial weak solutions with $I_\lambda(u_\lambda)=c_\lambda\ge\alpha>I_\lambda(v_\lambda)>0$; in particular, they are distinct.

\end{proof}

\begin{remark}
The result concerns a left neighborhood of $\lambda^*$ and does not establish multiplicity throughout $(\widehat\lambda^*,\lambda^*)$. The positivity here refers to the energy; we draw no conclusion about the pointwise sign of the solutions.
\end{remark}

\subsection*{Nonexistence of Nontrivial Solutions}
The existence and multiplicity of solutions for parameters in $(\widehat{\lambda}^*, \lambda^*-\varepsilon]$ remain unresolved. However, for $\lambda<\widehat{\lambda}^*$, we prove that the origin is the only weak solution. 

\begin{proof}[Proof of Theorem 1.2 \rm (iv)]
Let $\lambda<\widehat{\lambda}^*$ and suppose, by contradiction, that $u\in X\setminus\{0\}$ is a weak solution of \eqref{eq:problem}. Testing with $u$ gives $\langle I_\lambda'(u),u\rangle=0$, that is, $\lambda=\widehat{\lambda}_u(1)$, since $C(u)>0$. By the definition of $\widehat{\lambda}^*$,
\[
\widehat{\lambda}^*\le\inf_{t>0}\widehat{\lambda}_u(t)\le\widehat{\lambda}_u(1)=\lambda,
\]
which contradicts $\lambda<\widehat{\lambda}^*$. Therefore, \eqref{eq:problem} admits only the trivial weak solution for $\lambda<\widehat{\lambda}^*$.

\end{proof}

\section{Asymptotic Behavior of Solutions}\label{se.5}
We begin by recalling notation and definitions used throughout this section. Let \(I_0:X\to\mathbb{R}\) be the energy functional associated with problem \eqref{eq:problem} for \(\lambda=0\), given by
\[
I_0(u)=\frac{1}{r}\|u\|_X^r-\frac12 B(u).
\]
Also, we fix \(e\in E\setminus\{0\}\) such that
\[
I_0(e)<0.
\]
This condition is used to construct the admissible paths that characterize the mountain-pass levels. We
also recall that \(c_\lambda:=\inf_{\gamma\in\Gamma}\max_{t\in[0,1]}I_\lambda(\gamma(t))\), where
\(\Gamma:=\{\gamma\in C([0,1],X):\gamma(0)=0,\ \gamma(1)=e\}\) and \(e\) is the element fixed above.

\medskip

For every \(\lambda\in\mathbb{R}\), we define the Nehari set by
\[
\mathcal N_\lambda
:=
\left\{
u\in X\setminus\{0\}:
\langle I_\lambda'(u),u\rangle=0
\right\}.
\]
The positive and negative parts of \(\mathcal N_\lambda\) are defined, respectively, by
\[
\mathcal N_\lambda^+
:=
\left\{
u\in\mathcal N_\lambda:
\left.\frac{d}{dt}\langle I_\lambda'(tu),tu\rangle\right|_{t=1}>0
\right\},
\quad 
\mathcal N_\lambda^-
:=
\left\{
u\in\mathcal N_\lambda:
\left.\frac{d}{dt}\langle I_\lambda'(tu),tu\rangle\right|_{t=1}<0
\right\}.
\]
We also define the set of nonnegative mountain-pass solutions by
\[
\mathcal M_\lambda
:=
\left\{
u\in E\setminus\{0\}:
I_\lambda'(u)=0,\quad I_\lambda(u)=c_\lambda
\right\}.
\]
First, we show that, for \(\lambda<0\) sufficiently close to zero,
\[
c_\lambda
=\min_{u\in\mathcal N_\lambda^-}I_\lambda(u),
\qquad
\mathcal M_\lambda\cap\mathcal N_\lambda^-\neq\varnothing.
\]

We also show that the mountain-pass solutions can be chosen so that their norms remain uniformly bounded as \(\lambda\to0^-\). A similar uniform bound holds for the nonnegative mountain-pass solutions obtained in the cooperative regime as \(\lambda\to0^+\). These uniform bounds are the basis of the asymptotic analysis in both regimes.
%In addition, we will show that the mountain-pass solutions can be chosen so that their norms remain uniformly bounded as \(\lambda\to0^-\). This characterization will be fundamental to the study of the asymptotic behavior of the solutions.

\begin{lemma}[A bounded selection of mountain-pass solutions]\label{lem:bounded-mountain-pass-selection}
Assume $(H_1)-(H_3)$. Then, there exist $\varepsilon>0$ and $R>0$ such that, for every $\lambda\in(-\varepsilon,0)$, the mountain-pass solutions can be selected so that
\[
u_\lambda\in\mathcal M_\lambda\cap\mathcal N_\lambda^-, \qquad \|u_\lambda\|_X<R.
\]
Also,
\[
c_\lambda\to c_0>0
\qquad\text{as }\lambda\to0^-.
\]
\end{lemma}

\begin{proof}
Recall that there exist constants $S_B, S_C > 0$ such that
\[
B(u)\le S_B\|u\|_X^2,\qquad
0 \le C(u)\le S_C\|u\|_X^{2q} \quad \text{ for all } u \in X,
\]
in addition to the usual compactness properties of the functionals $B$ and $C$. 
To simplify the notation, we define the positive constants
\[
\alpha=\frac1r-\frac12,\qquad
\beta=\frac12-\frac1{2q},\qquad
\kappa=\alpha-\beta\frac{2-r}{2q-2}
=\frac{(2-r)(2q-r)}{4rq}>0.
\]
On $\mathcal{N}_\lambda^-$, where $\langle I_\lambda'(u),u\rangle=0$, the homogeneity condition yields 
\[
\frac{d}{dt}\langle I_\lambda'(tu),tu\rangle\Big|_{t=1}
=-(2-r)\|u\|_X^r-(2q-2)\lambda C(u).
\]
Consequently, for all $u\in \mathcal{N}_\lambda^-$, we obtain 
\begin{equation}\label{eq:nehari-selection-energy-bound}
I_\lambda(u)=\alpha \|u\|_X^r+\lambda\beta C(u)
\ge\kappa\|u\|_X^r,
\end{equation}
besides knowing that $\|u\|_X\ge\rho:=S_B^{-1/(2-r)}$. Direct differentiation yields
\begin{equation}\label{eq:nehari-selection-fiber-derivative}
\frac{d}{dt}\left(\frac{\langle I_\lambda'(tv),tv\rangle}{t^2}\right)
=t^{r-3}\left[-(2-r)\|v\|_X^r
 -(2q-2)\lambda C(v)t^{2q-r}\right].
\end{equation}
In particular, if $v\in\mathcal{N}_\lambda^-$, then $\langle I_\lambda'(v),v\rangle=0$. Therefore, the expression in \eqref{eq:nehari-selection-fiber-derivative} becomes strictly negative for $0<t\le1$, which allows us to conclude that 
\begin{equation}\label{eq:nehari-selection-fiber-increasing}
I_\lambda(tv)\le I_\lambda(v) \quad \text{ for all } t \in [0,1].
\end{equation}

Fix the constants $M$ and $R$ satisfying
\[
M>1+\max_{0\le t\le1}I_0(te),\qquad
R>\max\{\|e\|_X,\rho\},\qquad \kappa R^r>M+1.
\]
Thus, for all $\lambda<0$ with $|\lambda|$ sufficiently small, the following estimates are simultaneously valid:
\begin{gather}
I_\lambda(e)<0,\quad \langle I_\lambda'(e),e\rangle<0,\quad
\max_{0\le t\le1}I_\lambda(te)<M,\label{eq:nehari-selection-initial-choices}\\
-(2q-2)\lambda S_CR^{2q-r}<2-r,\qquad
-\lambda\beta S_CR^{2q}<1.\label{eq:nehari-selection-smallness}
\end{gather}
The first line uses the fact that $B(e)>2\|e\|_X^r/r>\|e\|_X^r$. The second line ensures that every nonzero point on the Nehari manifold, whose norm is bounded above by $R$, belongs to $\mathcal{N}_\lambda^-$. Since $\langle I_\lambda'(te),te\rangle>0$ in a neighborhood of $t=0$ and $\langle I_\lambda'(e),e\rangle<0$, there exists $t_\lambda \in (0,1)$ such that $t_\lambda e$ belongs to this branch. Therefore, we can define 
\[
m_\lambda:=\inf_{\mathcal{N}_\lambda^-}I_\lambda,
\]
which satisfies the inequalities $0<\kappa\rho^r\le m_\lambda<M$ and $c_\lambda<M$.

Let $(z_n)\subset\mathcal{N}_\lambda^-$ now be a minimizing sequence for the level $m_\lambda$. By \eqref{eq:nehari-selection-energy-bound}, for $n$ sufficiently large, we have the bound $\|z_n\|_X^r<M/\kappa<R^r$. Passing to a subsequence, we obtain 
\[
z_n\rightharpoonup z,\qquad B(z_n)\to B(z),\qquad C(z_n)\to C(z),
\qquad \|z_n\|_X^r\to a:=B(z)+\lambda C(z).
\]
By the lower semicontinuity of the norm, we know that $\|z\|_X^r\le a$, which implies $\langle I_\lambda'(z),z\rangle\le0$. In addition, since $B(z)\ge\rho^r$, we conclude that the limit $z\ne0$. Thus, there exists some $t\in(0,1]$ such that $tz\in\mathcal{N}_\lambda^-$. Suppose, by contradiction, that the strict inequality $\|z\|_X^r<a$ holds. Then, relying on \eqref{eq:nehari-selection-fiber-increasing}, we would have
\[
I_\lambda(tz)
=\lim_{n \to \infty} I_\lambda(tz_n)-\frac{t^r}{r}(a-\|z\|_X^r)
<\lim_{n \to \infty} I_\lambda(tz_n)\le m_\lambda,
\]
which is a contradiction to the definition of the infimum $m_\lambda$. Thus, we must have $\|z\|_X^r=a$, so $z_n\to z$ strongly in $X$. Therefore, $z$ is a minimizer that achieves $m_\lambda$.

Using $B(|z|)\ge B(z)$ and $C(|z|)=C(z)$, we choose $s\in(0,1]$ such that $s|z|\in\mathcal{N}_\lambda^-$. Thus, 
\[
m_\lambda\le I_\lambda(s|z|)\le I_\lambda(sz)\le I_\lambda(z)=m_\lambda.
\]
From this we deduce the existence of a nonnegative minimizer $u_\lambda \in \mathcal{N}_\lambda^-$. Denoting the derivative functional by $\phi_\lambda(u) = \langle I_\lambda'(u), u \rangle$, the Lagrange multipliers rule gives us $I_\lambda'(u_\lambda)=\eta \phi_\lambda'(u_\lambda)$ for some multiplier $\eta \in \mathbb{R}$. Testing this equation with $u_\lambda$ and using that $\langle \phi_\lambda'(u_\lambda),u_\lambda\rangle<0$ because $u_\lambda\in\mathcal{N}_\lambda^-$, we obtain $\eta=0$. Therefore, 
\begin{equation}\label{eq:nehari-selection-critical-minimizer}
u_\lambda\in E\cap\mathcal{N}_\lambda^-,
\qquad I_\lambda'(u_\lambda)=0,\qquad I_\lambda(u_\lambda)=m_\lambda,
\qquad \rho\le\|u_\lambda\|_X<R.
\end{equation}

To show that $m_\lambda\le c_\lambda$, let us consider the following open set
\[
D_\lambda=\{0\}\cup\{u:\|u\|_X<R,\ \langle I_\lambda'(u),u\rangle>0\}.
\]
It is a neighborhood of the origin that does not contain $e$. On its boundary, the points with norm $\|u\|_X=R$ satisfy $\langle I_\lambda'(u),u\rangle\ge0$. By condition \eqref{eq:nehari-selection-smallness}, we have
\[
I_\lambda(u)=\alpha R^r+\lambda\beta C(u)+\tfrac12\langle I_\lambda'(u),u\rangle
\ge\alpha R^r+\lambda\beta S_CR^{2q}>M.
\]
Therefore, any continuous path in the set $\Gamma$ whose maximum energy remains strictly less than $M$ must intersect the surface $\langle I_\lambda'(u),u\rangle=0$ at some interior point where $0<\|u\|_X<R$, so it crosses $\mathcal{N}_\lambda^-$. Since $c_\lambda<M$, this establishes that $m_\lambda\le c_\lambda$.

To prove the reverse inequality, we define $w_{\lambda,s}=(1-s)u_\lambda+se$, for $0\le s\le1$. Since the term $B$ is nonnegative on the positive cone $E$, this implies that
\[
B(w_{\lambda,s})\ge(1-s)^2B(u_\lambda)+s^2B(e)
\ge b_*:=\tfrac12\min\{\rho^r,B(e)\}>0,
\]
and $\|w_{\lambda,s}\|_X\le R$. Now, we can choose a scalar constant $L>1$ sufficiently large so that $L^rR^r/r-L^2b_*/2<-1$. Since $C(Lw_{\lambda,s})\le S_C(LR)^{2q}$, then, for any $\lambda<0$ with $|\lambda|$ sufficiently small, we have
\[
I_\lambda(Lw_{\lambda,s})<0 \quad \text{ for all } s \in [0,1].
\]
In parallel, we observe that $I_0(te)\le t^rI_0(e)\le I_0(e)<0$ for all $1\le t\le L$, which implies that $I_\lambda(te)<0$ along this same interval whenever $\lambda$ is sufficiently close to zero. Imposing the final restriction
\[
-(2q-2)\lambda S_CR^{2q-r}L^{2q-r}<2-r,
\]
and noting that $C(u_\lambda)/\|u_\lambda\|_X^r\le S_CR^{2q-r}$, the relation expressed in \eqref{eq:nehari-selection-fiber-derivative} guarantees that the function $\langle I_\lambda'(tu),tu\rangle/t^2$ is strictly decreasing on the interval $(0,L]$, reaching zero only at $t=1$. Therefore, 
\[
\max_{0\le t\le L}I_\lambda(tu_\lambda)=I_\lambda(u_\lambda)=m_\lambda.
\]
The continuous path formed by the concatenation $0\to Lu_\lambda\to Le\to e$ therefore belongs to the set of paths $\Gamma$ (after an appropriate reparametrization), and reaches its maximum energy exactly at the level $m_\lambda$. This proves that $c_\lambda\le m_\lambda$, which gives the equality of the levels
\[
c_\lambda=m_\lambda,\qquad
u_\lambda\in\mathcal{M}_\lambda\cap\mathcal{N}_\lambda^-,
\qquad
\kappa\|u_\lambda\|_X^r\le c_\lambda<M.
\]
All the smallness conditions required in the proof depend only on constants fixed earlier, so they can hold simultaneously for all $\lambda\in(-\varepsilon, 0)$ once $\varepsilon>0$ is chosen suitably.

It remains to analyze the asymptotic behavior of the energies. By the  geometry of the functional $I_0$, the classical Mountain-Pass condition $c_0>0$ holds. The inequality $I_\lambda\ge I_0$ implies that $c_\lambda\ge c_0$. Also, for each $\delta>0$, if we choose a fixed path $\gamma_\delta\in\Gamma$ such that $\max_t I_0(\gamma_\delta(t))\le c_0+\delta$, we obtain 
\[
c_\lambda\le c_0+\delta -
\frac{\lambda}{2q}\max_t C(\gamma_\delta(t)).
\]
This last maximum is finite, since $C$ is continuous and the image of the path is compact. By first taking the limit as $\lambda\to0^-$ and subsequently letting $\delta\to0^+$, we conclude that $c_\lambda\to c_0$.

\end{proof}

\begin{lemma}\label{lem:mountain_pass_nehari} 
Assume $(H_1)$ and $(H_2)$. Suppose also that there exists $e\in E\setminus\{0\}$ such that $I_0(e)<0$. Then,
\[
c_0=\inf_{u\in\mathcal{N}_0}I_0(u).
\]
\end{lemma}

\begin{proof}
Recall that there exists a constant $S_B>0$ such that $|B(u)|\le S_B\|u\|_X^2$. Choose \(\rho\in(0,\|e\|_X)\) sufficiently small so that \(S_B\rho^{2-r}<1\). Then
\[
\langle I_0'(u),u\rangle
\geq \rho^r-S_B\rho^2>0
\qquad\text{if }\|u\|_X=\rho.
\]
Since \(I_0(e)<0\), we also have
\[
B(e)>\frac2r\|e\|_X^r>\|e\|_X^r,
\]
and hence \(\langle I_0'(e),e\rangle<0\). Fix \(\gamma\in\Gamma\) and let
\[
t_1:=\max\{t\in[0,1]:\|\gamma(t)\|_X=\rho\}.
\]
Then \(\|\gamma(t)\|_X>\rho\) for \(t\in(t_1,1]\).
By continuity, the function $t\longmapsto\langle I_0'(\gamma(t)),\gamma(t)\rangle$ changes sign between \(t_1\) and \(1\). Thus there exists
\(t_0\in(t_1,1)\) such that
\[
\langle I_0'(\gamma(t_0)),\gamma(t_0)\rangle=0.
\]
Since \(\|\gamma(t_0)\|_X>\rho\), we obtain \(\gamma(t_0)\in\mathcal{N}_0\). Consequently,
\[
\max_{t\in[0,1]}I_0(\gamma(t))
\geq \inf_{u\in\mathcal{N}_0}I_0(u).
\]
Taking the infimum over \(\gamma\in\Gamma\) yields
\[
c_0\geq\inf_{u\in\mathcal{N}_0}I_0(u).
\]
For any \(u\in\mathcal{N}_0\), we have \(B(u)=\|u\|_X^r>0\). In addition, assumptions \((H_1)-(H_2)\) imply that $B(|u|)\geq B(u)>0$. Set
\[
t_u:=\left(\frac{B(u)}{B(|u|)}\right)^{1/(2-r)}
\in(0,1].
\]
Then
\[
\langle I_0'(t_u|u|),t_u|u|\rangle
=t_u^rB(u)-t_u^2B(|u|)=0,
\]
so \(t_u|u|\in\mathcal{N}_0\cap E\). Also, for every
\(z\in\mathcal{N}_0\),
\[
I_0(z)=\left(\frac1r-\frac12\right)\|z\|_X^r.
\]
Since \(t_u\leq1\), it follows that
\(I_0(t_u|u|)\leq I_0(u)\). Therefore,
\[
\inf_{u\in\mathcal{N}_0}I_0(u)
=\inf_{u\in\mathcal{N}_0\cap E}I_0(u).
\]
Fix \(w\in\mathcal{N}_0\cap E\). Since
\(B(w)=\|w\|_X^r\), we have
\[
I_0(tw)=\left(\frac{t^r}{r}-\frac{t^2}{2}\right)\|w\|_X^r, \qquad
\frac{d}{dt}I_0(tw)
=\|w\|_X^r(t^{r-1}-t).
\]
Thus \(I_0(tw)\) attains its global maximum at \(t=1\), that is, 
\[
\max_{t\geq0}I_0(tw)=I_0(w)>0.
\]
For \(s\in[0,1]\), define \(z_s=(1-s)w+se\). Since \(w,e\in E\setminus\{0\}\), we have \(B(z_s)>0\). By continuity and compactness, there exist \(\beta,M>0\) such that
\[
B(z_s)\geq\beta,
\qquad
\|z_s\|_X\leq M
\quad\text{for every }s\in[0,1].
\]
Consequently, choosing \(L>1\) sufficiently large, we obtain
\[
I_0(Lz_s)
\leq \frac{L^r}{r}M^r-\frac{L^2}{2}\beta<0
\qquad\text{for every }s\in[0,1],
\]
since \(r<2\). In addition, since \(B(e)>\|e\|_X^r\), for every
\(t\geq1\) we have
\[
\frac{d}{dt}I_0(te)
=t^{r-1}\|e\|_X^r-tB(e)
\leq t\bigl(\|e\|_X^r-B(e)\bigr)<0.
\]
Hence,
\[
I_0(te)\leq I_0(e)<0
\qquad\text{for every }t\geq1.
\]

The concatenation of the three segments $0\to Lw\to Le\to e$ defines, after a suitable reparametrization, a continuous path \(\gamma\in\Gamma\). Along the first segment, the energy attains its maximum at \(w\), since
\[
\max_{0\leq t\leq L} I_0(tw)=I_0(w).
\]
On the other hand, the energy is strictly negative along the remaining two segments, by the preceding estimates. Since \(I_0(w)>0\), it follows that
\[
\max_{t\in[0,1]}I_0(\gamma(t))=I_0(w).
\]
Therefore, by the definition of the mountain-pass level, $c_0\leq I_0(w)$. Taking the infimum over \(w\in\mathcal{N}_0\cap E\), we conclude that
\[
c_0
\leq\inf_{w\in\mathcal{N}_0\cap E}I_0(w)
=\inf_{u\in\mathcal{N}_0}I_0(u).
\]
Combining these two inequalities, we obtain $c_0=\inf_{u\in\mathcal{N}_0}I_0(u)$. This completes the proof.

\end{proof}

%\begin{proposition}\label{prop:asymptotic-behavior-negative}
%Assume $(H_1)-(H_3)$. For each $\lambda\in(\lambda^*,0)$, let $v_\lambda\in X$ be a global minimizer of $I_\lambda$ such that $I_\lambda(v_\lambda)<0$. Also, for $\lambda<0$ sufficiently close to zero, let $u_\lambda$ be the nonnegative mountain-pass solution chosen so that
%\[
%u_\lambda\in\mathcal M_\lambda\cap\mathcal N_\lambda^-,
%\]
%as provided by Lemma \ref{lem:bounded-mountain-pass-selection}. Then,
%\[
%\|v_\lambda\|_X\to+\infty \qquad\text{as }\lambda\to0^-.
%\]
%In addition, for every sequence $\lambda_n\to0^-$, there exist a subsequence, still denoted by $(\lambda_n)$, and a function $u_0\in E\setminus\{0\}$ such that
%\[
%u_{\lambda_n}\to u_0 \qquad\text{strongly in }X.
%\]
%This limit $u_0$ satisfies
%\[
%I_0'(u_0)=0 \quad \text{and} \quad I_0(u_0)=c_0.
%\]
%In particular, $u_0$ is a nontrivial, nonnegative weak solution of \eqref{eq:problem} at $\lambda=0$.
%\end{proposition}

\begin{proof}[Proof of Theorem \ref{thm:asymptotic-behavior-zero}]
Set
\[
\alpha:=\frac1r-\frac12>0,\qquad \beta:=\frac12-\frac1{2q}>0,\qquad
s_0:=\Big(\frac{c_0}{\alpha}\Big)^{1/r}.
\]
Recall that $|B(u)|\le S_B\|u\|_X^2$ and $0\le C(u)\le S_C\|u\|_X^{2q}$ for all $u\in X$. For every critical point $u \in X$ of $I_\lambda$, we have
\begin{equation}\label{eq:energy-identity}
I_\lambda(u)=I_\lambda(u)-\frac12\langle I_\lambda'(u),u\rangle
=\alpha\|u\|_X^r+\lambda\beta\,C(u).
\end{equation}

\medskip\noindent\emph{Proof of {\rm(i)}.}
Let $\lambda\in(\lambda^*,0)$. Since $v_\lambda$ is a global minimizer of the $C^1$ functional $I_\lambda$, we have $I_\lambda'(v_\lambda)=0$ and, by Proposition \ref{prop:minimizer}, $I_\lambda(v_\lambda)=\inf_XI_\lambda<0$; in particular $v_\lambda\ne0$. Then \eqref{eq:energy-identity} and the estimate $C(v_\lambda)\le S_C\|v_\lambda\|_X^{2q}$ give
\[
\alpha\|v_\lambda\|_X^r<(-\lambda)\beta\,C(v_\lambda)\le(-\lambda)\beta S_C\|v_\lambda\|_X^{2q}.
\]
As $v_\lambda\ne0$, dividing by $\|v_\lambda\|_X^r$ yields
\[
\|v_\lambda\|_X>\left(\frac{\alpha}{(-\lambda)\beta S_C}\right)^{1/(2q-r)}\to+\infty
\qquad\text{as }\lambda\to0^-.
\]

\medskip\noindent\emph{Proof of {\rm(ii)}.}
Let $e\in E\setminus\{0\}$ be the fixed element with $I_0(e)<0$, used as the endpoint of the paths defining $c_\lambda$ and $d_\lambda$, and let $\varepsilon,R>0$ be given by Lemma \ref{lem:bounded-mountain-pass-selection}. Fix $\lambda_1>0$ and let $\lambda_n\to0^-$ or $\lambda_n\to0^+$. Discarding finitely many indices, we may assume $\lambda_n\in(-\varepsilon,0)$ or $\lambda_n\in(0,\lambda_1]$, respectively. Write $u_n:=u_{\lambda_n}$, chosen as in Lemma \ref{lem:bounded-mountain-pass-selection} if $\lambda_n<0$, and given by Theorem \ref{teo:1.1} if $\lambda_n>0$ (in which case it is a weak solution in $X$, hence a critical point of $I_{\lambda_n}$ on $X$). Let $b_n:=I_{\lambda_n}(u_n)$ be its energy level, that is, $b_n=c_{\lambda_n}$ if $\lambda_n<0$ and $b_n=d_{\lambda_n}$ if $\lambda_n>0$.

If $\lambda_n<0$, then $\|u_n\|_X<R$ by Lemma~\ref{lem:bounded-mountain-pass-selection}. If $\lambda_n>0$, then $C(u)\ge0$ for all $u \in X$ gives $I_{\lambda_n}(u)\le I_0(u)$, so
$\max_tI_{\lambda_n}(\gamma(t))\le\max_tI_0(\gamma(t))$ for every path $\gamma \in \Gamma$. Taking the infimum over the paths and using the identification of the minimax levels in the proof of Theorem \ref{teo:1.1} (valid for every $\lambda\ge0$, including $\lambda=0$), we get $d_{\lambda_n}\le d_0=c_0$. Since $I_{\lambda_n}'(u_n)=0$, identity~\eqref{eq:energy-identity} gives
\begin{equation}\label{1.3-1}
\alpha\|u_n\|_X^r=d_{\lambda_n}-\lambda_n\beta\,C(u_n)\le d_{\lambda_n}\le c_0,
\end{equation}
that is, $\|u_n\|_X\le s_0$. Arguing as in the verification of the mountain-pass geometry in the proof of Theorem \ref{teo:1.1}, we choose $\rho\in(0,\|e\|_X)$ such that
\[
\frac{S_B}{2}\rho^{2-r}+\frac{\lambda_1S_C}{2q}\rho^{2q-r}\le\frac1{2r}.
\]
Then $I_\lambda(u)\ge\rho^r/(2r)$ whenever $\|u\|_X=\rho$ and $0\le\lambda\le\lambda_1$, and since every path $\gamma \in \Gamma$ crosses the sphere $\|u\|_X=\rho$, we obtain
\[
d_{\lambda_n}\ge\delta_0:=\frac{\rho^r}{2r}>0\qquad\text{for all }n.
\]

In both cases $(u_n)$ is bounded in $X$. Passing to a subsequence,
\[
u_n\rightharpoonup u_0\ \text{ in }X,\qquad
u_n\to u_0\ \text{ in }L^s(\Omega)\quad(1\le s<r^*).
\]
Since $I_{\lambda_n}'(u_n)=0$ and $I_{\lambda_n}=I_0-\frac{\lambda_n}{2q}C$, we have $I_0'(u_n)=\frac{\lambda_n}{2q}C'(u_n)$. The map $C'$ is bounded on bounded subsets of $X$ and $\lambda_n\to0$, hence
\begin{equation}\label{eq:nehari-selection-limit-derivative}
I_0'(u_n)\to0\qquad\text{in }X^*.
\end{equation}
Since $\langle I_0'(u_n), u_n-u_0 \rangle \to 0$, and by Lemma \ref{lem:convergences}, we have
\[
\int_\Omega|\nabla u_n|^{r-2}\nabla u_n\nabla(u_n-u_0)\,dx\to0.
\]
The $(S_+)$ property of the $r$-Laplacian gives $u_n\to u_0$ strongly in $X$. By the continuity of $I_0'$ and \eqref{eq:nehari-selection-limit-derivative}, we obtain $I_0'(u_0)=0$. Since $E$ is closed and $u_n\in E$, also $u_0\in E$.

Now, since $(u_n)$ is bounded, $C(u_n)$ is also bounded. Consequently, 
\[
I_0(u_n)=b_n+\frac{\lambda_n}{2q}C(u_n),\qquad \frac{\lambda_n}{2q}C(u_n)\to0.
\]
As $u_n\to u_0$ strongly in $X$, this gives $b_n\to I_0(u_0)$. If $\lambda_n<0$, then $b_n=c_{\lambda_n}\to c_0$ by Lemma \ref{lem:bounded-mountain-pass-selection}, so $I_0(u_0)=c_0>0$ and $u_0\ne0$. If $\lambda_n>0$, then $\delta_0\le b_n\le c_0$ gives $\delta_0\le I_0(u_0)\le c_0$, hence $u_0\ne0$. Since $I_0'(u_0)=0$, we have $u_0\in\mathcal N_0$, and Lemma \ref{lem:mountain_pass_nehari} yields $c_0=\inf_{\mathcal N_0}I_0\le I_0(u_0)$. Therefore $I_0(u_0)=c_0$ in both cases, and $u_0\in E\setminus\{0\}$ is a nontrivial nonnegative weak solution of \eqref{eq:problem} at $\lambda=0$.

It remains to prove that $d_\lambda\to c_0$ as $\lambda\to0^+$ (for $\lambda\to0^-$ this is part of Lemma~\ref{lem:bounded-mountain-pass-selection}). Take any sequence $\lambda_n\to0^+$, so that $(d_{\lambda_n})$ lies in $[\delta_0,c_0]$. Any subsequence of $(d_{\lambda_n})$ comes from parameters that still tend to $0^+$ and stay in $(0,\lambda_1]$. The argument above, applied to mountain-pass solutions at these levels, yields a further subsequence whose levels converge to $I_0(u_0)=c_0$. So every subsequence of $(d_{\lambda_n})$ has a further subsequence converging to $c_0$, and therefore $d_{\lambda_n}\to c_0$. Since $(\lambda_n)$ was arbitrary, we conclude that \(d_\lambda\to c_0\).

For the norms, recall that $I_0'(u_0)=0$ and $I_0(u_0)=c_0$. The identity \eqref{eq:energy-identity} with $\lambda=0$ gives $c_0=\alpha\|u_0\|_X^r$, that is, $\|u_0\|_X=s_0$. The limit $u_0$ may depend on the sequence, but its norm does not. By the previous paragraphs, every sequence $\lambda_n\to0^{\pm}$ has a subsequence along which $u_{\lambda_n}\to u_0$ strongly in $X$, so $\|u_{\lambda_n}\|_X\to\|u_0\|_X=s_0$ along it. Now suppose $\lim_{\lambda\to0^{\pm}}\|u_\lambda\|_X=s_0$ failed. There would be $\eta>0$ and $\lambda_n\to0^{\pm}$ with $\big|\|u_{\lambda_n}\|_X-s_0\big|\ge\eta$ for all $n$, and then no subsequence of $(\|u_{\lambda_n}\|_X)$ could converge to $s_0$, which is a contradiction.

\end{proof}

\bibliographystyle{ieeetr}
\bibliography{ref}

\end{document}